\documentclass[11pt]{article}
\usepackage{amsfonts,latexsym,amsmath,amscd,geometry}
\usepackage{upgreek}
\usepackage{amssymb}
\usepackage{latexsym}
\usepackage{color}

\makeatletter \@addtoreset{equation}{section} \makeatother

\newcommand \nc{\newcommand}
\newtheorem{theorem}{Theorem}[section]
\newtheorem{lemma}[theorem]{Lemma}

\newtheorem{definition}[theorem]{Definition}

\nc{\ba}{\begin{array}}\nc{\ea}{\end{array}}
\nc{\be}{\begin{eqnarray}}\nc{\ee}{\end{eqnarray}}
\nc{\beq}{\begin{equation}}\nc{\eeq}{\end{equation}}
\nc{\bex}{\begin{eqnarray*}}\nc{\eex}{\end{eqnarray*}}
\nc{\btm}{\begin{theorem}} \nc{\etm}{\end{theorem}}
\nc{\blm}{\begin{lemma}} \nc{\elm}{\end{lemma}}
\nc{\R}{\mathbb{R}} 
\nc{\va}{\varphi}
\nc{\ve}{\varepsilon}
\def\x{\mathbf{x}}

\def\pf{\noindent{\bf Proof.\quad}}\def\endpf{\hfill$\Box$}
\def\di{\mbox{div\,}}
\def\curl{\mbox{curl\,}}

 \newcommand{\n}{{\bf n}}
\newcommand{\bu}{{\bf u}}
\newcommand{\y}{{\bf y}}
\newcommand{\p}{{\bf P}}
\newcommand{\h}{{\mathcal H}}
\newcommand{\rs}{\langle\hat\rho\rangle_\sigma}
\newcommand{\rds}{\langle\hat\rho_\delta\rangle_\sigma}

\begin{document}
\title{Existence of global weak solutions to the compressible Ericksen-Leslie system in dimension one}

\author{Huajun Gong\footnote{College of Mathematics and Statistics, Shenzhen University, Shenzhen 518060, Guangdong, China.} \quad Tao Huang\footnote{Department of Mathematics, Wayne State University, Detroit, MI 48202, USA.}  \quad Changyou Wang\footnote{Department of Mathematics, Purdue University, West Lafayette, IN 47907, USA.}\\
}
\date{}
\maketitle

\begin{abstract}
We consider the compressible Ericksen-Leslie system of liquid crystal flows in one dimension. A global weak solution is constructed with initial density $\rho_0\geq 0$ and $\rho_0\in L^\gamma$ for $\gamma>1$.

\end{abstract}

\section{Introduction}

Nematic liquid crystals are composed of rod-like molecules characterized by average alignment of  the long axes of neighboring molecules, which have simplest structures among various types of liquid crystals. The dynamic theory of nematic liquid crystals has been first proposed by Ericksen \cite{ericksen62} and Leslie \cite{leslie68} in the 1960's, which is a macroscopic continuum description of the time evolution of both flow velocity field and orientation order parameter of rod-like liquid crystals.

In this paper, we will study the compressible Ericksen-Leslie system of liquid crystal flows (see \cite{Morro09}, \cite{Anna-Liu} for modeling). Let $\Omega\subset \R^3$ be a bounded domain with smooth boundary, and $\mathbb S^2$ be the unit sphere in $\R^3$. The compressible Ericksen-Leslie system is given as follows
\begin{equation}\label{comlce}
\begin{cases}
\rho_t+\nabla\cdot(\rho \bu)=0,\\ 
\rho\dot \bu+\nabla P=\nabla\cdot\sigma-\nabla\cdot\left(\frac{\partial W}{\partial\nabla \n}\otimes\nabla \n\right),
\\
{\bf g}+\frac{\partial W}{\partial  \n}-\nabla\cdot\left(\frac{\partial W}{\partial\nabla \n}\right)=\lambda\n.
\end{cases}
\end{equation}
Here, $\rho(\x,t):\Omega\times(0,\infty)\rightarrow \R$ is the density, $\bu(\x,t):\Omega\times(0,\infty)\rightarrow \R^3$ is the fluid velocity field, $\n(\x,t):\Omega\times(0,\infty)\rightarrow \mathbb S^2$ is the orientation order parameters of nematic material. $\lambda$ is the Lagrangian multiplier of the constraint $|\n|=1$,  $\dot f=f_t+\bu\cdot\nabla f$ is the material derivative of function $f$, and $\mathbf{a}\otimes \mathbf{b}=\mathbf{a}\, \mathbf{b}^T$ for column vectors $\mathbf{a}$ and $\mathbf{b}$ in $\mathbb{R}^3$.

The macrostructure of the crystals has been determined by the Oseen-Frank energy density (cf. \cite{oseen33,frank58}). One may take the Oseen-Frank energy density in the compressible case as 
\begin{align}\label{OFE}\begin{split}
2W(\rho, \n,\nabla \n)=&\frac{2}{\gamma-1}\rho^{\gamma}+K_1(\di \n)^2+K_2(\n\cdot\curl \n)^2+K_3|\n\times\curl \n|^2\\
&+(K_2+K_4)[\mbox{tr}(\nabla \n)^2-(\di \n)^2 ],
\end{split}
\end{align} 
where $\gamma>1$, and $K_j$, $j=1,2,3$, are the positive constants representing splay, twist, and bend effects respectively, with
$K_2\geq |K_4|$, $2K_1\geq K_2+K_4$. Then the pressure can be given by the Maxwell relation 
$$P(\rho)=\rho W_{\rho}(\rho, \n,\nabla \n)-W(\rho, \n,\nabla \n).$$  
For simplicity, we only consider the case $K_1=K_2=K_3=1$, $K_4=0$ in this paper. The Oseen-Frank energy in the compressible case becomes 
$$
2W(\rho, \n,\nabla \n)=\frac{2}{\gamma-1}\rho^{\gamma}+|\nabla \n|^2.
$$ 
Therefore
$$
\nabla\cdot\left(\frac{\partial W}{\partial\nabla \n}\otimes\nabla \n\right)=\nabla\cdot\left(\nabla\n\odot\nabla\n\right),\quad \frac{\partial W}{\partial  \n}=0,\quad \nabla\cdot\left(\frac{\partial W}{\partial\nabla \n}\right)=\Delta \n,\quad P=\rho^{\gamma}-\frac12|\nabla n|^2.
$$
Let
$$
D= \frac12(\nabla \bu+\nabla^{T} \bu),\quad \omega= \frac12(\nabla \bu-\nabla^{T}\bu)=\frac{1}{2}\left(
\frac{\partial u^i}{\partial x_j}-\frac{\partial u^j}{\partial x_i}\right),\quad N=\dot \n-\omega \n, 
$$
represent the rate of strain tensor, skew-symmetric part of the strain rate, and the
rigid rotation part of director changing rate by fluid vorticity, respectively. The kinematic transport ${\bf g}$ is given by
\begin{align}\label{g}
{\bf g}=\gamma_1 N +\gamma_2D\n-\gamma_2(\n^TD\n)\n 
\end{align}
which represents the effect of the macroscopic flow field on the microscopic structure. The material coefficients $\gamma_1$ and $\gamma_2$ reflect the molecular shape and the slippery part between fluid and particles. The first term of ${\bf g}$ represents the rigid rotation of molecules, while the second term stands for the stretching of molecules by the flow.
The viscous (Leslie) stress tensor $\sigma$ has the following form
(cf. \cite{Les} \cite{Anna-Liu})
\begin{align}\label{sigma}\begin{split}
\sigma=& \alpha_0(\n^TD\n)\mathbb I+ \alpha_1 (\n^TD\n)\n\otimes \n +\alpha_2N\otimes \n+ \alpha_3 \n\otimes N\\
&   + \alpha_4D + \alpha_5(D\n)\otimes \n+\alpha_6\n\otimes (D\n)
+\alpha_7(\mbox{tr}\,D)\,\mathbb I+\alpha_8(\mbox{tr}\,D)\,\n\otimes\n.
\end{split}
\end{align}
These coefficients $\alpha_j$ $(0 \leq j \leq 8)$, depending on material and temperature, are called Leslie coefficients. The following relations are often assumed in the literature.
\begin{align}\label{a2g}
\gamma_1 =\alpha_3-\alpha_2,\quad \gamma_2 =\alpha_6 -\alpha_5,\quad \alpha_2+ \alpha_3 =\alpha_6-\alpha_5. 
\end{align}
The first two relations are compatibility conditions, while the third relation is called Parodi's relation, derived from Onsager reciprocal relations expressing the equality of certain relations between flows and forces in thermodynamic systems out of equilibrium (cf. \cite{Parodi70}). 
They also satisfy the following empirical relations (cf. \cite{Les}, \cite{Anna-Liu})
\begin{align}\label{alphas}
&\alpha_4>0,\quad 2\alpha_1+3\alpha_4+2\alpha_5+2\alpha_6>0,\quad \gamma_1=\alpha_3-\alpha_2>0,\\
&  2\alpha_4+\alpha_5+\alpha_6>0,\quad 4\gamma_1(2\alpha_4+\alpha_5+\alpha_6)>(\alpha_2+\alpha_3+\gamma_2)^2\notag\\
&\alpha_4+\alpha_7>\alpha_1+\frac{\gamma_2^2}{\gamma_1}\geq 0,\quad \notag\\
&2\alpha_4+\alpha_5+\alpha_6-\frac{\gamma_2^2}{\gamma_1}>\alpha_0+\alpha_1+\alpha_5+\alpha_6+\alpha_8\geq 0.\notag
\end{align}

It is easy to see that an example of coefficients $\alpha_1, \cdots, \alpha_8$ satisfying
\eqref{a2g} and \eqref{alphas}  can be taken as follows
$$
\alpha_0=\alpha_1=\alpha_5=\alpha_6=\alpha_7=\alpha_8=0,\quad \alpha_2=-1,\quad  \alpha_3= \alpha_4=1,
$$
so that 
$$
\gamma_1=\alpha_3-\alpha_2=2>0,\quad \gamma_2=\alpha_6-\alpha_5=\alpha_2+\alpha_3=0.
$$

A simplified compressible Ericksen-Leslie system has been recently studied. The idea of simplification was first proposed for the incompressible system by Lin in\cite{Lin89}. In dimension one, the global strong and weak solutions have been constructed in \cite{dlww12} and \cite{dww11}. In dimension two, under the assumption that the initial data of $\n$ is contained in $\mathbb S^2_+$, global weak solutions have been constructed in \cite{jsw13}. In dimension three, the local existence of strong solutions has been studied by \cite{hww121} and \cite{hww122}, and when the initial  data of $\n$ is contained in $\mathbb S^2_+$, global weak solutions have been constructed in \cite{llw15}. The incompressible limit of compressible nematic liquid crystal flows has been studied by \cite{dhxwz13}. 

We also mention a related work \cite{jlt19}, in which the Ericksen–Leslie’s parabolic–hyperbolic liquid crystal model has been studied. 
For small initial data, they have shown the existence of global  solutions in dimension three.

\subsection{One dimensional model and statement of main results}

One of the main motivations of this paper is to 
investigate the impact of general Leslie stress tensors to the solutions of the compressible Ericksen-Leslie system with coefficients satisfying algebraic conditions \eqref{a2g} and \eqref{alphas} ensuring the energy dissipation property.
Because of the technical complexity of the Ericksen-Leslie system in higher dimensions, 
we will only consider the following simpler case in one dimension, 
in which the director field $\n$ is assumed to map into  the equator $\mathbb S^1$,
$$
\bu=\big(u(x,t),\ v(x,t), 0\big)^T, \quad \n=\big(\cos n(x,t),\ \sin n(x,t),0 \big)^T
$$
for any $x\in [0,1]$ and $t\in (0,\infty)$. From the derivation given by Section 2 below, 
the system \eqref{comlce} becomes 
\begin{equation}\label{comlce1d}
\begin{cases}
\rho_t+(\rho u)_x=0,\\ 
(\rho u)_t+(\rho u ^2)_x+\big(\rho^{\gamma}\big)_x=J^1-n_{xx}n_x,
\\
(\rho v)_t+(\rho u v)_x=J^2,\\
\gamma_1\left(\dot n-\frac12v_x\right)
-\gamma_2\left(u_x\cos n\sin n+\frac12v_x(1-2\cos^2 n)\right) =n_{xx}.
\end{cases}
\end{equation}
Here 
\beq\notag
\begin{split}
J^1=&(\alpha_0+\alpha_5+\alpha_6+\alpha_8)\big(u_x\cos^2 n\big)_x+\alpha_1\big(u_x\cos^4n\big)_x-(\alpha_2+\alpha_3)\big(\dot n\cos n\sin n\big)_x+(\alpha_4+\alpha_7)u_{xx}\\
&+\alpha_0\big(v_x\cos n\sin n\big)_x+\alpha_1\big(v_x\cos^3n\sin n\big)_x+\frac12(\alpha_2+\alpha_3+\alpha_5+\alpha_6)\big(v_x\cos n\sin n\big)_x,
\end{split}
\eeq
and
\beq\notag
\begin{split}
J^2=&\alpha_1\big(u_x\cos^3n\sin n\big)_x+\alpha_2\big(\dot n\cos^2 n\big)_x-\alpha_3\big(\dot n\sin^2 n\big)_x+(\alpha_6+\alpha_8)\big(u_x\cos n\sin n\big)_x\\
&+\alpha_1\big(v_x\cos^2n\sin^2 n\big)_x+\frac12(-\alpha_2+\alpha_5)\big(v_x\cos^2 n\big)_x+\frac12(\alpha_3+\alpha_6)\big(v_x\sin^2 n\big)_x+\frac12\alpha_4v_{xx}.
\end{split}
\eeq
For this system, we consider the following initial and boundary values
\beq\label{1dinitial}
(\rho,\, \rho u,\, \rho v,\, n)(x,0)=(\rho_0,\, m_0,\, l_0,\, n_0)(x),
\eeq
\beq\label{1dbdyvalue}
u(0,t)=v(0,t)=u(1,t)=v(1,t)=0,\quad n_x(0,t)=n_x(1,t)=0.
\eeq

Denote the energy of the system \eqref{comlce1d} by
\beq\notag
\mathcal E(t):=\frac{1}{2}\int_0^1\rho (u^2+v^2)
+\frac{1}{\gamma-1}\int_0^1\rho^{\gamma}+\frac12\int_0^1n_x^2.
\eeq
For any smooth solution $(\rho, u, v, n)$, the energy functional satisfies the following energy inequality, whose proof will be provided in Section 3,
\beq\label{enest1}
\begin{split}
\frac{d}{dt}\mathcal E(t)
&=-\mathcal{D}\\
&:=-\int_0^1\left[\sqrt{\gamma_1}\dot n-\frac12\left(
 \frac{\gamma_2}{\sqrt{\gamma_1}}u_x\sin(2n)
 + \frac{1}{\sqrt{\gamma_1}}(\gamma_1-\gamma_2\cos(2n))v_x
\right)\right]^2\\
&-\int_0^1\left[\frac14\left(-\alpha_1-\frac{\gamma_2^2}{\gamma_1}\right)u_x^2+(\alpha_4+\alpha_7)u_{x}^2\right]
-\frac14\int_0^1\left(2\alpha_4+\alpha_5+\alpha_6-\frac{\gamma_2^2}{\gamma_1}\right)v_x^2\\
&-\frac14\left(\alpha_1+\frac{\gamma_2^2}{\gamma_1}\right)
\int_0^1\left(u_x\cos(2n)+v_x\sin(2n)\right)^2\\
&-(\alpha_0+\alpha_1+\alpha_5+\alpha_6+\alpha_8)\int_0^1\left[\big(u_x\cos n+\frac12v_x\sin n\big)^2-\frac14 v_x^2\sin^2 n\right].
\end{split}
\eeq
By the assumptions \eqref{alphas} on coefficients,   the system \eqref{comlce1d} is dissipative.

\begin{definition}\label{defweaksol}
For any time $0<T<\infty$, a collection of functions $(\rho, u, v, n)(x,t)$ is a global weak solution to the initial and boundary value problem \eqref{comlce1d}-\eqref{1dbdyvalue} if
\begin{itemize}
\item[(1)]
$$\rho\geq 0,\ \mbox {a.e.}, \quad \rho\in L^{\infty}(0,T; L^{\gamma}), \quad \rho u^2, \rho v^2 \in L^{\infty}(0,T; L^{1}), \quad  u, v \in L^{2}(0,T; H^{1}_0)$$
$$
n\in L^{\infty}(0,T; H^{1})\cap L^{2}(0,T;  H^{2}),\quad
n_t\in L^{2}(0,T;   L^{2}).
$$

\item[(2)] The equations of $\rho$, $u$, $v$ are satisfied in the weak sense, while the equation of $n$ is valid a.e.. The initial condition \eqref{1dinitial} is satisfied in the weak sense. 

\item[(3)]  The energy inequality is valid for a.e. $t\in(0,T)$
$$
\mathcal E(t)+\int_0^t\mathcal{D}\leq \mathcal E_0 =
\frac{1}{2}\int_0^1 \frac{m_0^2+l_0^2}{\rho_0}
+\frac{1}{\gamma-1}\int_0^1\rho_0^{\gamma}+\frac12\int_0^1(n_0)_x^2.
$$

\end{itemize}
\end{definition}

The following is the main results in this paper.

\begin{theorem}\label{mainth1} Assume that the coefficients of Leslie stress tensor satisfy the algebraic conditions \eqref{a2g} and \eqref{alphas}. Then, for any $0<T<\infty$ and any initial data
\beq\label{asmp-initial}
0\le \rho_0\in L^{\gamma}, \quad \frac{m_0}{\sqrt{\rho_0}}, \quad \frac{l_0}{\sqrt{\rho_0}}\in L^2, \quad n_0\in H^1,
\eeq
there is a global weak solution $(\rho, u, v, n)(x,t)$ on $(0,1)\times (0,T)$ to the initial and boundary value problem \eqref{comlce1d}-\eqref{1dbdyvalue}. Furthermore,
$\rho\in L^{2\gamma}((0,1)\times (0,T))$.
\end{theorem}

\medskip

The main ideas of the proof utilize and extend those from \cite{fnp01},  \cite{jiangzhang01}, and \cite{feireisl04} in the study of the compressible Navier-Stokes equations, where the quantity called {\it effective viscous flux} has played  crucial roles in controlling the oscillation of the density function $\rho$. However,  the general Leslie stress tensors
in the compressible Ericksen-Leslie system \eqref{comlce1d} induce two complicate second-order terms $J^1$ and $J^2$ that prohibit direct applications of the method of effective viscous flux. In this paper, we observe that
with the algebraic conditions \eqref{a2g} and \eqref{alphas},  the system of $\bu=(u,\, v)^T$ can still be shown
to be uniformly parabolic (see \eqref{postiveA} and \eqref{vec1dlce} below), i.e. the coefficient matrix of the second-order terms is uniformly elliptic. Using the inverse of coefficient matrix of the second-order terms, we can then define a modified form of effective viscous flux as in Lemma \ref{lemma5.3}, which yields the desired estimates
that are necessary in the limiting process of approximated solutions.

\medskip

The paper is organized as follows. In Section 2, we will sketch a derivation of the system \eqref{comlce1d}. In Section 3, we will derive some a priori estimates for smooth solutions of \eqref{comlce1d}. In Section 4, an approximated system will be introduced,
 and the existence of global regular solutions of this approximated system will be proven. In Section 5, we 
 will prove the existence of global weak solutions through some delicate analysis of the convergence process.

\section{Derivation of the model in one dimension}

This section is devoted to the derivation of the system \eqref{comlce1d} in dimension one. 
If a solution takes the form
$$
\bu=\big(u(x,t),\ v(x,t)\big)^T, \quad \n=\big(\cos n(x,t),\ \sin n(x,t) \big)^T, 
\ (x,t)\in (0,1)\times (0,T), 
$$
then
$$
\nabla \bu=
\left[
\begin{array}{cc}
u_x & 0\\
v_x &0
\end{array}
\right],
\quad
\nabla^T \bu=
\left[
\begin{array}{cc}
u_x & v_x\\
0 &0
\end{array}
\right],
$$
so that
$$
D=
\left[
\begin{array}{cc}
u_x & \frac12v_x\\
\frac12v_x &0
\end{array}
\right]
\quad 
\omega=
\left[
\begin{array}{cc}
0 & -\frac12v_x\\
\frac12v_x &0
\end{array}
\right],$$
$$
\mbox{tr}\,D=u_x,\quad  N=\dot \n-\omega\n=\left(\dot n-\frac12v_x\right) \big(-\sin n,\ \cos n \big)^T.
$$
Direct calculations imply that 
$$
D\n=\left(u_x\cos n+\frac12v_x\sin n,\ \frac12v_x\cos n\right)^T,\quad \n^TD\n=u_x\cos^2 n+v_x\cos n\sin n,
$$
$$
\n\otimes\n=\left[
\begin{array}{cc}
\cos^2 n& \cos n\sin n\\
\cos n\sin n&\sin^2 n
\end{array}
\right],
$$
$$
(\n^TD\n)\n\otimes \n=
(u_x\cos^2 n+v_x\cos n\sin n)
\left[
\begin{array}{cc}
\cos^2 n& \cos n\sin n\\
\cos n\sin n&\sin^2 n
\end{array}
\right],
$$
$$
N\otimes\n=\left(\dot n-\frac12v_x\right)
\left[
\begin{array}{cc}
-\cos n\sin n& -\sin^2 n\\
\cos^2 n&\cos n\sin n
\end{array}
\right],
$$
$$
\n\otimes N=\left(\dot n-\frac12v_x\right)
\left[
\begin{array}{cc}
-\cos n\sin n& \cos^2 n\\
-\sin^2 n&\cos n\sin n
\end{array}
\right],
$$
$$
(D\n)\otimes \n=\left[
\begin{array}{cc}
u_x\cos ^2n+\frac12v_x\cos n\sin n& u_x\cos n\sin n+\frac12v_x\sin^2 n\\
\frac12v_x\cos^2 n&\frac12v_x\cos n\sin n
\end{array}
\right],$$
$$
\n\otimes (D\n)=\left[
\begin{array}{cc}
u_x\cos ^2n+\frac12v_x\cos n\sin n&\frac12v_x\cos^2 n\\
 u_x\cos n\sin n+\frac12v_x\sin^2 n&\frac12v_x\cos n\sin n
\end{array}
\right].
$$
Hence 
$$
\nabla \cdot\sigma=\big(J^1, J^2\big)^T
$$
where 
\beq\notag
\begin{split}
J^1=&(\alpha_0+\alpha_5+\alpha_6+\alpha_8)\big(u_x\cos^2 n\big)_x+\alpha_1\big(u_x\cos^4n\big)_x-(\alpha_2+\alpha_3)\big(\dot n\cos n\sin n\big)_x+(\alpha_4+\alpha_7)u_{xx}\\
&+\alpha_0\big(v_x\cos n\sin n\big)_x+\alpha_1\big(v_x\cos^3n\sin n\big)_x+\frac12(\alpha_2+\alpha_3+\alpha_5+\alpha_6)\big(v_x\cos n\sin n\big)_x,
\end{split}
\eeq
and
\beq\notag
\begin{split}
J^2=&\alpha_1\big(u_x\cos^3n\sin n\big)_x+\alpha_2\big(\dot n\cos^2 n\big)_x-\alpha_3\big(\dot n\sin^2 n\big)_x+(\alpha_6+\alpha_8)\big(u_x\cos n\sin n\big)_x\\
&+\alpha_1\big(v_x\cos^2n\sin^2 n\big)_x+\frac12(-\alpha_2+\alpha_5)\big(v_x\cos^2 n\big)_x+\frac12(\alpha_3+\alpha_6)\big(v_x\sin^2 n\big)_x+\frac12\alpha_4v_{xx}.
\end{split}
\eeq
The terms related to $\n$ can be computed as follows
$$
\n_t=n_t \big(-\sin n,\ \cos n \big)^T,
$$
$$
\n_x=n_x \big(-\sin n,\ \cos n \big)^T,\quad |\n_x|^2=(n_x)^2
$$
$$
\bu\cdot\n=u\n_x=un_x \big(-\sin n,\ \cos n \big)^T
$$
$$
\n_{xx}=n_{xx} \big(-\sin n,\ \cos n \big)^T+(n_x)^2 \big(-\cos n,\ -\sin n \big)^T,
$$
\beq\notag
\begin{split}
\nabla\cdot\left(\nabla\n\odot\nabla\n\right)-\frac12\nabla|\nabla \n|^2
=&\Delta \n\nabla \n=\big(n_{xx}n_x,\ 0\big)^T.
\end{split}
\eeq
Therefore, $u(x,t)$ satisfies
\beq
\rho u_t+\rho u u_x+\big(\rho^{\gamma}\big)_x=J^1-n_{xx}n_x,
\eeq
and $v(x,t)$ satisfies
\beq
\rho v_t+\rho u v_x=J^2.
\eeq

Now we can calculate the equation of $n$ as follows.
\beq\notag
\begin{split}
{\bf g}=&\gamma_1 N +\gamma_2D\n-\gamma_2(\n^TD\n)\n \\
&=\gamma_1\left(\dot n-\frac12v_x\right) \big(-\sin n,\ \cos n \big)^T+\gamma_2\left(u_x\cos n+\frac12v_x\sin n,\ \frac12v_x\cos n\right)^T\\
&-\gamma_2 \big(u_x\cos^2 n+v_x\cos n\sin n\big)\big(\cos n,\ \sin n\big)^T\\
&=\gamma_1\left(\dot n-\frac12v_x\right) \big(-\sin n,\ \cos n \big)^T\\
&+\gamma_2\left(u_x\cos n\sin^2 n+\frac12v_x\sin n(1-2\cos^2 n),\ -u_x\cos^2n\sin n+\frac12v_x\cos n(1-2\sin^2 n)\right)^T\\
&=\gamma_1\left(\dot n-\frac12v_x\right) \big(-\sin n,\ \cos n \big)^T\\
&-\gamma_2\left(u_x\cos n\sin n+\frac12v_x(1-2\cos^2 n)\right) \big(-\sin n,\ \cos n \big)^T,
\end{split}
\eeq
$$
\lambda \n=\left(|\nabla \n|^2+\gamma_1 N\cdot\n\right)\n=(n_x)^2 \big(\cos n,\ \sin n\big)^T.
$$
Therefore $n(x,t)$ satisfies
\beq
\gamma_1\left(\dot n-\frac12v_x\right)
-\gamma_2\left(u_x\cos n\sin n+\frac12v_x(1-2\cos^2 n)\right) =n_{xx}.
\eeq
Thus the system \eqref{comlce} reduces to \eqref{comlce1d}.


\section{A priori estimates}

In this section, we will prove several useful a priori estimates for smooth solutions of system \eqref{comlce1d}. 

\begin{lemma}\label{lemma1}
Any smooth solution to the system \eqref{comlce1d} satisfies the following energy inequality 
\beq\label{enest2}
\begin{split}
\frac{d}{dt}\mathcal E(t)
=&-\int_0^1\left[\sqrt{\gamma_1}\dot n-\frac12\left(
 \frac{\gamma_2}{\sqrt{\gamma_1}}u_x\sin(2n)
 + \frac{1}{\sqrt{\gamma_1}}(\gamma_1-\gamma_2\cos(2n))v_x
\right)\right]^2\\
&-\int_0^1\left[\frac14\left(-\alpha_1-\frac{\gamma_2^2}{\gamma_1}\right)u_x^2+(\alpha_4+\alpha_7)u_{x}^2\right]
-\frac14\int_0^1\left(2\alpha_4+\alpha_5+\alpha_6-\frac{\gamma_2^2}{\gamma_1}\right)v_x^2\\
&-\frac14\left(\alpha_1+\frac{\gamma_2^2}{\gamma_1}\right)
\int_0^1\left(u_x\cos(2n)+v_x\sin(2n)\right)^2\\
&-(\alpha_0+\alpha_1+\alpha_5+\alpha_6+\alpha_8)\int_0^1\left[\big(u_x\cos n+\frac12v_x\sin n\big)^2-\frac14 v_x^2\sin^2 n\right].
\end{split}
\eeq
\end{lemma}

\pf
Multiplying the second equation by $u$, the third equation by $v$ and integrating over $[0,1]$, we have 
\beq\notag
\begin{split}
\frac{1}{2}\frac{d}{dt}\int_0^1\rho (u^2+v^2)
+\frac{1}{\gamma-1}\frac{d}{dt}\int_0^1\rho^{\gamma}
=\int_0^1\left(J^1u+J^2v-un_{xx}n_x\right).
\end{split}
\eeq
Multiplying the last equation by $\dot n$ and integrating over $[0,1]$, we obtain
\beq\notag
\begin{split}
\frac{d}{dt}\frac12\int_0^1(n_x)^2
+\gamma_1\int_0^1\dot n^2
=\int_0^1\left[\frac12\gamma_2u_x\sin(2n)\dot n+\frac12(\gamma_1-\gamma_2\cos(2n))v_x\dot n+un_{xx}n_x\right].
\end{split}
\eeq
Adding these two equations together, we have 
\beq\label{pfsec3.1}
\begin{split}
&\frac{1}{2}\frac{d}{dt}\int_0^1\rho (u^2+v^2)
+\frac{1}{\gamma-1}\frac{d}{dt}\int_0^1\rho^{\gamma}
+\frac12\frac{d}{dt}\int_0^1(n_x)^2
\\
=&\int_0^1\left(J^1u+J^2v\right)-\gamma_1\int_0^1\dot n^2+\int_0^1\frac12\left[\gamma_2u_x\sin(2n)\dot n+(\gamma_1-\gamma_2\cos(2n))v_x\dot n\right].
\end{split}
\eeq
By integrating by parts, we can estimate the term related to $J^1, J^2$ as follows
\beq\label{pfsec3.2}
\begin{split}
&\int_0^1J^1u\\
=&-\int_0^1\left[(\alpha_0+\alpha_5+\alpha_6+\alpha_8)u_x^2\cos^2 n+\alpha_1u_x^2\cos^4n+(\alpha_4+\alpha_7)u_{x}^2\right]\\
&-\int_0^1\left[\alpha_0u_xv_x\cos n\sin n+\alpha_1u_xv_x\cos^3n\sin n+\frac12(\alpha_2+\alpha_3+\alpha_5+\alpha_6)u_xv_x\cos n\sin n\right]\\
&+\int_0^1(\alpha_2+\alpha_3)u_x\dot n\cos n\sin n,
\end{split}
\eeq
\beq\label{pfsec3.3}
\begin{split}
&\int_0^1J^2v\\
=&-\int_0^1\left[\alpha_1v_x^2\cos^2n\sin^2 n+\frac12(-\alpha_2+\alpha_5)v_x^2\cos^2 n+\frac12(\alpha_3+\alpha_6)v_x^2\sin^2 n+\frac12\alpha_4v_{x}^2\right]\\
&-\int_0^1\left[\alpha_1u_xv_x\cos^3n\sin n+(\alpha_6+\alpha_8)u_xv_x\cos n\sin n\right]\\
&-\int_0^1\left[\alpha_2v_x\dot n\cos^2 n-\alpha_3v_x\dot n\sin^2 n\right].
\end{split}
\eeq
First notice that all the terms related to $\alpha_1$ in \eqref{pfsec3.2} and \eqref{pfsec3.3} can be written as 
\beq\label{pfsec3.4}
\begin{split}
-\alpha_1\int_0^1\left[u_x^2\cos^4n+2u_xv_x\cos^3n\sin n+v_x^2\cos^2n\sin^2 n\right]\\
=-\alpha_1\int_0^1\left[u_x\cos^2n+v_x\cos n\sin n\right]^2.
\end{split}
\eeq
The other term related to $u_xv_x$ in \eqref{pfsec3.2} and \eqref{pfsec3.3} (without terms with $\alpha_1$) can be written as
\beq\label{pfsec3.5}
\begin{split}
&-\int_0^1\left[\alpha_0u_xv_x\cos n\sin n+\frac12(\alpha_2+\alpha_3+\alpha_5+\alpha_6)u_xv_x\cos n\sin n+(\alpha_6+\alpha_8)u_xv_x\cos n\sin n\right]\\
=&-\int_0^1u_xv_x\cos n\sin n\left[\alpha_0+\frac12(\alpha_2+\alpha_3+\alpha_5+\alpha_6)+(\alpha_6+\alpha_8)\right]\\
=&-\int_0^1\left(\alpha_0+2\alpha_6+\alpha_8\right)u_xv_x\cos n\sin n,
\end{split}
\eeq
where we have used $\alpha_2+\alpha_3=\alpha_6-\alpha_5$.
The terms related to $u_x^2$, $v_x^2$ in \eqref{pfsec3.2} and \eqref{pfsec3.3} (without terms with $\alpha_1$) can be written as 
\beq\label{pfsec3.6}
\begin{split}
&-\int_0^1\left[(\alpha_0+\alpha_5+\alpha_6+\alpha_8)u_x^2\cos^2 n+(\alpha_4+\alpha_7)u_{x}^2\right]\\
&-\int_0^1\left[\frac14(2\alpha_4-\alpha_2+\alpha_5+\alpha_3+\alpha_6)v_x^2-\frac12\gamma_2v_x^2\cos (2n)\right].
\\
\end{split}
\eeq
What left in \eqref{pfsec3.1}-\eqref{pfsec3.3} are all terms related to $u_x\dot n$ and $v_x\dot n$
\beq\label{pfsec3.7}
\begin{split}
&\int_0^1\left[\frac12\gamma_2u_x\sin(2n)\dot n+(\alpha_2+\alpha_3)u_x\dot n\cos n\sin n\right]\\
&+\int_0^1\left[\frac12(\gamma_1-\gamma_2\cos(2n))v_x\dot n-\alpha_2v_x\dot n\cos^2 n+\alpha_3v_x\dot n\sin^2 n\right]\\
=&\int_0^1\gamma_2u_x\dot n\sin(2n)
+\int_0^1(\gamma_1-\gamma_2\cos(2n))v_x\dot n,
\end{split}
\eeq
where we have used $\gamma_1=\alpha_3-\alpha_2$ and $\gamma_2=\alpha_2+\alpha_3=\alpha_6-\alpha_5$. 
Therefore, putting \eqref{pfsec3.4}-\eqref{pfsec3.7} into \eqref{pfsec3.1}, we obtain
\beq\label{pfsec3.8}
\begin{split}
&\frac{1}{2}\frac{d}{dt}\int_0^1\rho (u^2+v^2)
+\frac{1}{\gamma-1}\frac{d}{dt}\int_0^1\rho^{\gamma}+\frac{d}{dt}\frac12\int_0^1(n_x)^2
\\
=&-\alpha_1\int_0^1\left[u_x\cos^2n+v_x\cos n\sin n\right]^2-\int_0^1u_xv_x\cos n\sin n\left(\alpha_0+2\alpha_6+\alpha_8\right)
\\
&-\int_0^1\left[(\alpha_0+\alpha_5+\alpha_6+\alpha_8)u_x^2\cos^2 n+(\alpha_4+\alpha_7)u_{x}^2\right]\\
&-\int_0^1\left[\frac14(2\alpha_4+\alpha_5+\alpha_6+\gamma_1)v_x^2-\frac12\gamma_2v_x^2\cos (2n)\right]\\
&-\gamma_1\int_0^1\dot n^2+\int_0^1\gamma_2u_x\dot n\sin(2n)
+\int_0^1(\gamma_1-\gamma_2\cos(2n))v_x\dot n.
\end{split}
\eeq
We first complete the square for all terms with $\dot n$ in \eqref{pfsec3.4} 
\beq\label{pfsec3.9}
\begin{split}
&\gamma_1\int_0^1\dot n^2-\int_0^1\gamma_2u_x\dot n\sin(2n)
-\int_0^1(\gamma_1-\gamma_2\cos(2n))v_x\dot n\\
=&\gamma_1\int_0^1\dot n^2-2\cdot\frac12\int_0^1\sqrt{\gamma_1}\dot n
\left(
 \frac{\gamma_2}{\sqrt{\gamma_1}}u_x\sin(2n)
 + \frac{1}{\sqrt{\gamma_1}}(\gamma_1-\gamma_2\cos(2n))v_x
\right)\\
=&\int_0^1\left[\sqrt{\gamma_1}\dot n-\frac12\left(
 \frac{\gamma_2}{\sqrt{\gamma_1}}u_x\sin(2n)
 + \frac{1}{\sqrt{\gamma_1}}(\gamma_1-\gamma_2\cos(2n))v_x
\right)\right]^2\\
&-\frac{1}{4}\int_0^1 \left(
 \frac{\gamma_2}{\sqrt{\gamma_1}}u_x\sin(2n)
 + \frac{1}{\sqrt{\gamma_1}}(\gamma_1-\gamma_2\cos(2n))v_x
\right)^2.
\end{split}
\eeq
The last term in \eqref{pfsec3.9} can also be rewritten as follows
\beq\label{pfsec3.10}
\begin{split}
 &\left(
 \frac{\gamma_2}{\sqrt{\gamma_1}}u_x\sin(2n)
 + \frac{1}{\sqrt{\gamma_1}}(\gamma_1-\gamma_2\cos(2n))v_x
\right)^2\\
=&\frac{\gamma_2^2}{\gamma_1}u_x^2\sin^2(2n)
+2\frac{\gamma_2}{\gamma_1}u_xv_x\sin(2n)(\gamma_1-\gamma_2\cos(2n))
+ \frac{1}{\gamma_1}(\gamma_1-\gamma_2\cos(2n))^2v_x^2\\
=&\frac{\gamma_2^2}{\gamma_1}u_x^2\sin^2(2n)
+2u_xv_x\sin(2n)\left(\gamma_2-\frac{\gamma_2^2}{\gamma_1}\cos(2n)\right)\\
&+ \left(\gamma_1-2\gamma_2\cos(2n)+\frac{\gamma_2^2}{\gamma_1}\cos^2(2n)\right)v_x^2.
\end{split}
\eeq

To complete the square for the remaining terms, we first investigate the terms 
containing $u_xv_x$ in \eqref{pfsec3.9} and \eqref{pfsec3.10}:
\beq\label{pfsec3.11}
\begin{split}
&\frac12\alpha_1\int_0^1u_xv_x\sin(2n)(1+\cos(2n))+\frac12\int_0^1\left(\alpha_0+2\alpha_6+\alpha_8\right)u_xv_x\sin (2n)\\
&-\frac{1}{2}\int_0^1 u_xv_x\sin(2n)\left(\gamma_2-\frac{\gamma_2^2}{\gamma_1}\cos(2n)\right)\\
=&\frac12\int_0^1\left(\alpha_0+\alpha_1+\alpha_5+\alpha_6+\alpha_8\right)u_xv_x\sin (2n)
+\frac12\int_0^1\left(\alpha_1+\frac{\gamma_2^2}{\gamma_1}\right)u_xv_x\sin(2n)\cos(2n)\\
=&\int_0^1\left(\alpha_0+\alpha_1+\alpha_5+\alpha_6+\alpha_8\right)u_xv_x\sin n\cos n
+\frac12\int_0^1\left(\alpha_1+\frac{\gamma_2^2}{\gamma_1}\right)u_xv_x\sin(2n)\cos(2n).
\end{split}
\eeq
Thus we can calculate the terms containing $u_x^2$ in \eqref{pfsec3.9} and \eqref{pfsec3.10} as follows
\beq\label{pfsec3.12}
\begin{split}
&\frac14\int_0^1\left[\alpha_1u_x^2(1+\cos(2n))^2-\frac{\gamma_2^2}{\gamma_1}u_x^2\sin^2(2n)\right]
\\
&+\int_0^1\left[(\alpha_0+\alpha_5+\alpha_6+\alpha_8)u_x^2\cos^2 n+(\alpha_4+\alpha_7)u_{x}^2\right]\\
=&\frac14\int_0^1\left[\alpha_1u_x^2(1+2\cos(2n)+\cos^2(2n))-\frac{\gamma_2^2}{\gamma_1}u_x^2+\frac{\gamma_2^2}{\gamma_1}u_x^2\cos^2(2n)\right]
\\
&+\int_0^1\left[(\alpha_0+\alpha_5+\alpha_6+\alpha_8)u_x^2\cos^2n +2(\alpha_4+\alpha_7)u_{x}^2\right]\\
=&\frac14\int_0^1\left(\alpha_1+\frac{\gamma_2^2}{\gamma_1}\right)u_x^2\cos^2(2n)+\int_0^1(\alpha_0+\alpha_1+\alpha_5+\alpha_6+\alpha_8)u_x^2\cos^2n
\\
&+\int_0^1\left[\frac14\left(-\alpha_1-\frac{\gamma_2^2}{\gamma_1}\right)u_x^2+(\alpha_4+\alpha_7)u_{x}^2\right].
\end{split}
\eeq
Similarly, the terms involving $v_x^2$ in \eqref{pfsec3.9} and \eqref{pfsec3.10} can be calculated as follows
\beq\label{pfsec3.13}
\begin{split}
&\frac14\int_0^1\alpha_1v_x^2\sin^2(2n)+\int_0^1\left[\frac14(2\alpha_4+\alpha_5+\alpha_6+\gamma_1)v_x^2-\frac12\gamma_2v_x^2\cos (2n)\right]\\
&-\frac14\int_0^1 \left(\gamma_1-2\gamma_2\cos(2n)+\frac{\gamma_2^2}{\gamma_1}\cos^2(2n)\right)v_x^2\\
=&\frac14\int_0^1\alpha_1v_x^2\sin^2(2n)+\frac14\int_0^1\left(2\alpha_4+\alpha_5+\alpha_6-\frac{\gamma_2^2}{\gamma_1}\cos^2(2n)\right)v_x^2\\
=&\frac18\int_0^1\big(2\alpha_1+3\alpha_4+2\alpha_5+2\alpha_6\big)v_x^2\sin^2(2n)+\frac18\int_0^1\alpha_4v_x^2\sin^2(2n)\\
&+\frac14\int_0^1\left(2\alpha_4+\alpha_5+\alpha_6-\frac{\gamma_2^2}{\gamma_1}\right)v_x^2\cos^2(2n).
\end{split}
\eeq
For the terms with coefficient $\alpha_1+\frac{\gamma_2^2}{\gamma_1}$ in \eqref{pfsec3.11} and \eqref{pfsec3.12}, we have
\beq\label{pfsec3.14}
\begin{split}
&\frac14\int_0^1\left(\alpha_1+\frac{\gamma_2^2}{\gamma_1}\right)u_x^2\cos^2(2n)+\frac12\int_0^1\left(\alpha_1+\frac{\gamma_2^2}{\gamma_1}\right)u_xv_x\sin(2n)\cos(2n)\\
=&\frac14\left(\alpha_1+\frac{\gamma_2^2}{\gamma_1}\right)
\int_0^1\big[\left(u_x\cos(2n)+v_x\sin(2n)\right)^2-v_x^2\sin^2(2n)\big].
\end{split}
\eeq
The terms with coefficient $\alpha_0+\alpha_1+\alpha_5+\alpha_6+\alpha_8$ \eqref{pfsec3.11} and \eqref{pfsec3.12} can be written as 
\beq\label{pfsec3.15}
\begin{split}
&(\alpha_0+\alpha_1+\alpha_5+\alpha_6+\alpha_8)\int_0^1\big(u_x^2\cos^2n+u_xv_x\sin n\cos n\big)\\
=&(\alpha_0+\alpha_1+\alpha_5+\alpha_6+\alpha_8)\int_0^1\left[\big(u_x\cos n+\frac12v_x\sin n\big)^2-\frac14 v_x^2\sin^2 n\right].
\end{split}
\eeq
Collecting all the terms involving $v_x^2$ in \eqref{pfsec3.13}-\eqref{pfsec3.15}, we have
\beq\label{pfsec3.16}
\begin{split}
&\frac18\int_0^1\big(2\alpha_1+3\alpha_4+2\alpha_5+2\alpha_6\big)v_x^2\sin^2(2n)+\frac18\int_0^1\alpha_4v_x^2\sin^2(2n)\\
&+\frac14\int_0^1\left(2\alpha_4+\alpha_5+\alpha_6-\frac{\gamma_2^2}{\gamma_1}\right)v_x^2\cos^2(2n)
-\frac14\left(\alpha_1+\frac{\gamma_2^2}{\gamma_1}\right)\int_0^1v_x^2\sin^2(2n)\\
&-\frac14(\alpha_0+\alpha_1+\alpha_5+\alpha_6+\alpha_8)\int_0^1v_x^2\sin^2 n\\
=&\frac14\int_0^1\left(2\alpha_4+\alpha_5+\alpha_6-\frac{\gamma_2^2}{\gamma_1}\right)v_x^2-\frac14(\alpha_0+\alpha_1+\alpha_5+\alpha_6+\alpha_8)\int_0^1v_x^2\sin^2 n.
\end{split}
\eeq
Therefore, putting the identities \eqref{pfsec3.9} \eqref{pfsec3.14}-\eqref{pfsec3.16} into \eqref{pfsec3.8} 
yields 
\beq\notag
\begin{split}
&\frac{1}{2}\frac{d}{dt}\int_0^1\rho (u^2+v^2)
+\frac{1}{\gamma-1}\frac{d}{dt}\int_0^1\rho^{\gamma}+\frac{d}{dt}\frac12\int_0^1(n_x)^2
\\
=&-\int_0^1\left[\sqrt{\gamma_1}\dot n-\frac12\left(
 \frac{\gamma_2}{\sqrt{\gamma_1}}u_x\sin(2n)
 + \frac{1}{\sqrt{\gamma_1}}(\gamma_1-\gamma_2\cos(2n))v_x
\right)\right]^2\\
&-\int_0^1\left[\frac14\left(-\alpha_1-\frac{\gamma_2^2}{\gamma_1}\right)u_x^2+(\alpha_4+\alpha_7)u_{x}^2\right]
-\frac14\int_0^1\left(2\alpha_4+\alpha_5+\alpha_6-\frac{\gamma_2^2}{\gamma_1}\right)v_x^2\\
&-\frac14\left(\alpha_1+\frac{\gamma_2^2}{\gamma_1}\right)
\int_0^1\left(u_x\cos(2n)+v_x\sin(2n)\right)^2\\
&-(\alpha_0+\alpha_1+\alpha_5+\alpha_6+\alpha_8)\int_0^1\left[\big(u_x\cos n+\frac12v_x\sin n\big)^2-\frac14 v_x^2\sin^2 n\right],
\end{split}
\eeq
which completes the proof of Lemma. 
\endpf

From the energy inequality above, we can obtain the following estimates for $n$. 
\begin{lemma}\label{lemma2}
For any smooth solution to the system \eqref{comlce1d}, it holds that
\beq\label{estnxxnt}
\|n_{xx}\|_{L^2(0,T;L^2)}+\| n_t\|_{L^2(0,T;L^2)}\leq C(\mathcal E_0,T).
\eeq
\end{lemma}

\pf First notice that the equation of $n$ is 
\beq
\gamma_1\left(\dot n-\frac12v_x\right)
-\gamma_2\left(u_x\cos n\sin n+\frac12v_x(1-2\cos^2 n)\right) =n_{xx}.
\eeq
It is not hard to see that 
\beq\notag
\begin{split}
&\gamma_1\left(\dot n-\frac12v_x\right)
-\gamma_2\left(u_x\cos n\sin n+\frac12v_x(1-2\cos^2 n)\right)\\
=&\gamma_1\dot n-\frac12\gamma_2u_x\sin(2n)-\frac12(\gamma_1-\gamma_2\cos(2n))v_x.
\end{split}
\eeq
By the energy inequality, we obtain the estimates for $n_{xx}$. Next,  by the equation of $n$ and 
the energy inequality, we obtain the estimate for $ n_t$.
\endpf

We also need to show the higher integrability of $\rho$, which is inspired by the argument in \cite{dww11}.

\begin{lemma}\label{lemma3}
For any smooth solution to the system \eqref{comlce1d}, it holds that
\beq\label{estrho2g}
\|\rho\|_{L^{2\gamma}([0,1]\times[0,T];)}\leq C(\mathcal E_0,T).
\eeq
\end{lemma}

\pf
First set
$$G(x,t):=\int_0^x\rho^\gamma-x\int_0^1\rho^\gamma.$$ 
It is easy to see that 
$$
\frac{\partial G}{\partial x}=\rho^\gamma-\int_0^1\rho^\gamma,\quad G(0,t)=G(1,t)=0.
$$
Notice that the equation of $u$ can be written as  
\beq\notag
(\rho u)_t+(\rho u^2)_x+\big(\rho^{\gamma}\big)_x=J^1-\frac{1}{2}((n_x)^2)_x
\eeq
where 
\beq\notag
\begin{split}
J^1=&(\alpha_0+\alpha_5+\alpha_6+\alpha_8)\big(u_x\cos^2 n\big)_x+\alpha_1\big(u_x\cos^4n\big)_x-(\alpha_2+\alpha_3)\big(\dot n\cos n\sin n\big)_x+(\alpha_4+\alpha_7)u_{xx}\\
&+\alpha_0\big(v_x\cos n\sin n\big)_x+\alpha_1\big(v_x\cos^3n\sin n\big)_x+\frac12(\alpha_2+\alpha_3+\alpha_5+\alpha_6)\big(v_x\cos n\sin n\big)_x.
\end{split}
\eeq
Multiplying this equation by $G(x,t)$, integrating over $[0,1]\times(0,T)$, and using integrating by parts, 
we obtain that
\beq\label{rho2-1}
\begin{split}
\int_0^T\int_0^1\rho^{2\gamma}=&\int_0^T\left(\int_0^1\rho^{\gamma}\right)^2+\int_0^T\int_0^1(\rho u)_tG(x,t)-\int_0^T\int_0^1 \rho u^2\frac{\partial G(x,t)}{\partial x}\\
&-\int_0^T\int_0^1J_1G(x,t)-\frac12\int_0^T\int_0^1|n_x|^2\frac{\partial G(x,t)}{\partial x}\\
=&\sum_{i=1}^5I_i.
\end{split}
\eeq
For the first term, it is easy to estimate by energy inequality 
\beq\notag
I_1\leq C(\mathcal E_0,T).
\eeq
For the second term, we need use integrating by parts with respect to $t$ to obtain  
\beq\notag
\begin{split}
I_2=&\int_0^1\rho uG(x,T)-\int_0^1\rho uG(x,0)-\int_0^T\int_0^1\rho uG_t(x,t)\\
&\leq C\sup\limits_{0\leq t\leq T}\left(\int_0^1\rho|u|\int_0^1\rho^\gamma\right)-\int_0^T\int_0^1\rho uG_t(x,t)\\
&\leq C\sup\limits_{0\leq t\leq T}\left(\int_0^1\rho|u|^2\int_0^1\rho^\gamma+\int_0^1\rho\int_0^1\rho^\gamma\right)-\int_0^T\int_0^1\rho uG_t(x,t)\\
&\leq C(\mathcal E_0,T)-\int_0^T\int_0^1\rho uG_t(x,t).
\end{split}
\eeq 
To estimate the last term here, we multiply the equation of $\rho$ by $\gamma\rho^{\gamma-1}$ to get
\beq\notag
(\rho^{\gamma})_t+(\rho^{\gamma} u)_x+(\gamma-1)\rho^{\gamma} u_x=0.
\eeq
Then it holds
\beq\notag
\begin{split}
&-\int_0^T\int_0^1\rho uG_t(x,t)\\
=&-\int_0^T\int_0^1\rho u\left(\int_0^x\rho^\gamma_t-x\int_0^1\rho^\gamma_t\right)\\
=&\int_0^T\int_0^1\rho u\int_0^x\left((\rho^{\gamma} u)_x+(\gamma-1)\rho^{\gamma} u_x\right)-\int_0^T\int_0^1x\rho u\int_0^1\left((\rho^{\gamma} u)_x+(\gamma-1)\rho^{\gamma} u_x\right)\\
=&\int_0^T\int_0^1\rho^{\gamma+1} u^2+(\gamma-1)\int_0^T\int_0^1\rho u\left(\int_0^x\rho^{\gamma} u_x-x\int_0^1\rho^{\gamma} u_x\right)\\
\leq &\int_0^T\int_0^1\rho^{\gamma+1} u^2+C\int_0^T\int_0^1\rho |u|\int_0^1\rho^{\gamma} |u_x|\\
\leq &\int_0^T\int_0^1\rho^{\gamma+1} u^2+C\int_0^T\left(\int_0^1(\rho+ \rho|u|^2)\left(\int_0^1\rho^{2\gamma}\right)^{\frac12} \left(\int_0^1|u_x|^2\right)^{\frac12}\right)\\
\leq &\int_0^T\int_0^1\rho^{\gamma+1} u^2+C(\mathcal E_0,T)\int_0^T\left(\left(\int_0^1\rho^{2\gamma}\right)^{\frac12} \left(\int_0^1|u_x|^2\right)^{\frac12}\right)\\
\leq &\int_0^T\int_0^1\rho^{\gamma+1} u^2+\frac14\int_0^T\int_0^1\rho^{2\gamma}+C(\mathcal E_0,T)\int_0^T\int_0^1|u_x|^2\\
\leq &\int_0^T\int_0^1\rho^{\gamma+1} u^2+\frac14\int_0^T\int_0^1\rho^{2\gamma}+C(\mathcal E_0,T),
\end{split} 
\eeq
where we have used the Cauchy inequality, the H$\ddot{\mbox{o}}$lder inequality, the Young inequality and the energy inequality.  
Hence we obtain
\beq\notag
\begin{split}
I_2\leq \int_0^T\int_0^1\rho^{\gamma+1} u^2+\frac14\int_0^T\int_0^1\rho^{2\gamma}+C(\mathcal E_0,T).
\end{split}
\eeq 
For the third term in \eqref{rho2-1}, it holds
\beq\notag
\begin{split}
I_3=-\int_0^T\int_0^1 \rho u^2\left(\rho^\gamma-\int_0^1\rho^\gamma\right)=-\int_0^T\int_0^1 \rho^{\gamma+1} u^2+C(\mathcal E_0,T).
\end{split}
\eeq
Then 
\beq\notag
I_2+I_3\leq \frac14\int_0^T\int_0^1\rho^{2\gamma}+C(\mathcal E_0,T).
\eeq
For the fourth term in \eqref{rho2-1}, by integration by parts  it holds
\beq\notag
\begin{split}
I_4\leq &\int_0^T\int_0^1 \left(|u_x|+|\dot{n}|+|v_x|\right)\rho^{\gamma}+\int_0^T\int_0^1\left(|u_x|+|n_t|+|v_x|\right)\int_0^1\rho^\gamma\\
\leq &\frac14\int_0^T\int_0^1\rho^{2\gamma} +C\int_0^T\int_0^1\left(|u_x|^2+|\dot{n}|^2+|v_x|^2\right)\\
&+C(\mathcal E_0,T)\int_0^T\int_0^1\left(|u_x|^2+|n_t|^2+|v_x|^2\right)+C(\mathcal E_0,T)\\
\leq &\frac14\int_0^T\int_0^1\rho^{2\gamma}+C(\mathcal E_0,T).
\end{split}
\eeq
For the last term in \eqref{rho2-1}, it holds
\beq\notag
\begin{split}
I_5= -\frac12\int_0^T\int_0^1 |n_x|^2\left(\rho^\gamma-\int_0^1\rho^\gamma\right)\leq C(\mathcal E_0,T).
\end{split}
\eeq
Therefore, by adding all the estimates together in \eqref{rho2-1} we obtain
 \beq\notag
\begin{split}
\int_0^T\int_0^1\rho^{2\gamma}
\leq\frac12\int_0^T\int_0^1\rho^{2\gamma} +C(\mathcal E_0,T),
\end{split}
\eeq
which implies the estimate \eqref{estrho2g}. \endpf


\section{Approximated solutions}
In this section, we first consider the case that the initial values are smooth enough, i.e.
$\rho_0\in C^{1}$, $u_0, v_0, n_0\in C^{2}$, and $0< c_0^{-1}\leq \rho_0\leq c_0$ and $u_0=\frac{m_0}{\rho_0}$, $v_0=\frac{l_0}{\rho_0}$, and then construct the Galerkin approximation of $\rho$, $u$, $v$ and $n$.

\bigskip
\noindent \textbf{Step 1}. Recall that 
$$
\phi_j(x)=\sin\left(j\pi x\right), \quad j=1,2,...
$$
is an orthogonal base of $L^2(0,1)$. For any positive integer $k$, set
$$
\mathcal X_k=\mbox{span}\{\phi_1,\,\phi_2,\,\cdots\, \phi_k\}.
$$
and
\beq\notag
u_0^k=\sum_{j=0}^k \bar c_j^k\phi_j(x),\quad v_0^k=\sum_{j=0}^k\bar d_j^k\phi_j(x),
\eeq  
for some constants 
$$
\bar c_j^k=\int_0^1u_0\phi_j,\quad \bar d_j^k=\int_0^1v_0\phi_j.
$$ 
Then $(u_0^k,\,v_0^k)\rightarrow (u_0,\,v_0)$ in $C^2$ as $k\rightarrow \infty$.
Let
\beq\notag
u_k=\sum_{j=0}^kc_j^k(t)\phi_j(x),\quad v_k=\sum_{j=0}^kd_j^k(t)\phi_j(x)
\eeq  
be the finite dimensional approximation of $u$, and $v$, and we want to solve the approximation system:

\begin{equation}\label{applce}
\begin{cases}
(\rho_k)_t+(\rho_k u_k)_x=0,\\ 
\rho_k (u_k)_t+\rho_k u_k (u_k)_x+\big(\rho_k^{\gamma}\big)_x=J^1_k-(n_k)_{xx}(n_k)_x,
\\
\rho_k (v_k)_t+\rho_k u_k (v_k)_x=J^2_k,\\
\gamma_1\left(\dot n_k-\frac12(v_k)_x\right)
-\gamma_2\left((u_k)_x\cos n_k\sin n_k+\frac12(v_k)_x(1-2\cos^2 n_k)\right) =(n_k)_{xx}.
\end{cases}
\end{equation}
Here $J_k^1$, $J_k^2$ have the same form as $J^1$, $J^2$, but with $u, v$ replaced by $u_k$, $v_k$.  
For this system, we consider the following initial and boundary values
\beq\label{appinitial}
(\rho_k,\, u_k,\,  v_k,\, n_k)(x,0)=(\rho_0,\, u^k_0,\, v^k_0,\, n_0)(x),
\eeq
\beq\label{appbdyvalue}
u_k(0,t)=v_k(0,t)=u_k(1,t)=v_k(1,t)=0,\quad (n_k)_x(0,t)=(n_k)_x(1,t)=0.
\eeq

\bigskip
\noindent \textbf{Step 2}. The first step is to solve $\rho_k$ and $n_k$ by assuming $u_k, v_k\in C^{0}(0,T; C^{2})$ for a fixed $k$. To this end, we rewrite the equations of $\rho_k$ and $n_k$ in the Lagrange coordinate system.

Without loss of generality, in this section, we assume that 
\beq\label{exrho_0}
\int_0^1\rho_0(x)\,dx=1.
\eeq
For any $T>0$, we introduce the Lagrangian coordinate $(X,\tau)\in (0,1)\times [0,T)$ by
\beq\notag
X(x,t)=\int_0^x\rho_k(y,t)\,dy, \quad \tau(x,t)=t.
\eeq
If $\rho_k(x,t)\in C^1((0,1)\times[0,T))$ is positive
and $\int_0^1\rho_k(x,t)\,dx=1$ for all $t\in [0,T)$, then the map $(x,t)\rightarrow (X,\tau):(0,1)\times (0,T)
\to (0,1)\times (0,T) $ is a $C^1$-bijection such that $X(0,t)=0,\ X(1,t)=1$. 
By the chain rule, we have
\beq\notag
\frac{\partial}{\partial t}=-\rho_k u_k\frac{\partial}{\partial X}+\frac{\partial}{\partial \tau},\quad \frac{\partial}{\partial x}=\rho_k\frac{\partial}{\partial X}.
\eeq
The equation of $\rho_k$ can be rewritten as 
\begin{equation}\label{applceL1}
(\rho_k)_{\tau}+\rho_k^2 (u_k)_X=0,
\end{equation}
along with the initial condition
\beq\label{Linitial1}
\rho_k(X,0)=\rho_0.
\eeq

 Suppose $u_k\in C^{0}(0,T;C^{2})$ with $\|u_k\|_{C^{0}(0,T;C^{2})}\leq M_0$.
Then $\rho_k$ can be solved explicitly by
\beq\label{rhoXT}
\rho_k(X,\tau)=\frac{\rho_0(X)}{1+\rho_0(X)\int_0^{\tau}(u_k)_X(X,s)\,ds}.
\eeq
Hence, for any $T\leq \frac{1}{2c_0M_0}$, we have
\beq\label{2c}
\rho_k(X,\tau)
\leq \frac{\rho_0(X)}{1-\left|\rho_0(X)\int_0^{\tau}(u_k)_X(X,s)\,ds\right|}
\leq \frac{c_0}{1-c_0M_0T}\leq 2c_0,
\eeq
\beq\label{2c1}
\rho_k(X,\tau)
\geq \frac{\rho_0(X)}{1+\left|\rho_0(X)\int_0^{\tau}(u_k)_X(X,s)\,ds\right|}
\geq \frac{c_0^{-1}}{1+c_0M_0T}\geq \frac{c_0^{-1}}{2}.
\eeq
Similarly, since $\rho_0\in C^{1}$, $u_k\in C^{0}(0,T; C^{2})$, we conclude that for sufficiently small
$T(c_0, M_0)>0$,
\beq\label{M1}
\|\rho_k\|_{C^{0}(0,T; C^{1})}+\| (\rho_k)_t\|_{C^{0}((0,1)\times(0,T))}\leq M_1,
\eeq 
for some positive constant $M_1$.

Furthermore, suppose that $\rho^1_k, \rho^2_k$ are solutions of equation \eqref{applceL1} 
corresponding to $u_k^1, u^2_k \in C^{0}(0,T; C^{2})$,
with the same initial condition, we can conclude from \eqref{applceL1} that 
\beq\notag
\left(\frac{1}{\rho_k^1}-\frac{1}{\rho_k^2}\right)_\tau=\left(u^1_k-u^2_k\right)_X.
\eeq
Integrating with respect to $\tau$, we obtain
\beq\notag
\rho^1_k-\rho^2_k=\rho^1_k\rho^2_k\int_0^\tau\left(u^1_k-u^2_k\right)_X
\eeq
which, combined with \eqref{M1},  implies that
\beq\label{Liprho}
\|\rho^1_k-\rho^2_k\|_{C^{0}(0,T; C^{1})}
\leq C(M_1, T)T\|u^1_k-u^2_k\|_{C^{0}(0,T; C^{2})}.
\eeq

\bigskip
\noindent \textbf{Step 3}. Similarly, we can rewrite the equation of $n$ in the Lagrange coordinate as

\begin{equation}\label{applceL2}
\gamma_1\left((n_k)_\tau-\frac12\rho_k(v_k)_X\right)
-\frac{\gamma_2} {2}\left(\rho_k(u_k)_X\sin (2n_k)-\rho_k(v_k)_X\cos (2n_k)\right) =\rho_k\big(\rho_k(n_k)_{X}\big)_X.
\end{equation}
For this system, we consider the following initial and boundary values 
\beq\label{Linitial2}
n_k(X,0)=n_0(X),
\eeq
\beq\label{Lbdyvalue2}
 (n_k)_X(0,\tau)=(n_k)_X(1,\tau)=0.
\eeq
By the standard Schauder theory of parabolic equations, we conclude that
\beq\label{M2}
\begin{split}
\|n_k\|_{C^{1}(0,T; C^{2})}
\leq&C\|n_0\|_{C^{2}}+C\|\rho_k (v_k)_X\|_{C^{0}((0,1)\times(0,T))}+C\|\rho_k (u_k)_X\|_{C^{0}((0,1)\times(0,T))}
\leq  M_2,
\end{split}
\eeq 
for some positive constant $M_2$.

Furthermore, suppose that $n^1_k, n^2_k$ are solutions of equation \eqref{applceL2} 
corresponding to $\rho^1_k, \rho^2_k \in C^{1}((0,1)\times(0,T))$ and $u_k^1, u^2_k \in C^{0}(0,T; C^{2})$,
subject to the same initial condition. Denote 
$$
\bar n_k=n^1_k-n^2_k,\quad \bar \rho_k=\rho^1_k-\rho^2_k,\quad \bar u_k=u^1_k-u^2_k.
$$
Then from \eqref{applceL2} we have that 
\beq\notag
\begin{split}
&\gamma_1 (\bar n_k)_\tau-(\rho_k^1)^2(\bar n_k)_{XX}\\
&=
\bar\rho_k(\rho^1_k+\rho^2_k)(n_k^2)_{XX}+\bar\rho_k(\rho_k^1)_X(n_k^1)_{X}+\rho_k^2(\bar\rho_k)_X(n_k^1)_{X}+\rho_k^2(\rho_k^2)_X(\bar n_k)_{X}\\
&+\frac{\gamma_1}{2}\left(\bar\rho_k(v_k^1)_X+\rho_k^2(\bar v_k^1)_X\right)\\
&-\frac{\gamma_2} {2}\left(\bar\rho_k(v_k^1)_X\cos (2n_k^1)+\rho_k^2(\bar v_k)_X\cos (2n_k^1)-2\rho_k^2(v_k^2)_X\sin (\bar n_k)\sin (n_k^1+n_k^2)\right) \\
&+\frac{\gamma_2} {2}\left(\bar\rho_k(u_k^1)_X\sin (2n_k^1)+\rho_k^2(\bar u_k)_X\sin (2n_k^1)+2\rho_k^2(u_k^2)_X\sin (\bar n_k)\cos (n_k^1+n_k^2)\right).
\end{split}
\eeq
By the standard $W^{2,1}_2$-estimate of parabolic equations, we conclude that
\beq\notag
\begin{split}
&\|\bar n_k\|_{W^{2,1}_2([0,1]\times (0,T))}\\
&\leq  C\|\bar \rho_k\|_{L^2(0,T; H^{1})}+C\|\bar n_k\|_{L^{2}(0,T; L^{2})}+C\|\bar v_k\|_{L^2(0,T; H^{1})}+C\|\bar u_k\|_{L^2(0,T; H^{1})}\\
&\leq CT^{\frac12}\|\bar \rho_k\|_{C^{0}(0,T; C^{1})}+C\|\bar n_k\|_{L^{2}(0,T; L^{2})}+CT^{\frac12}\|\bar v_k\|_{C^{0}(0,T; C^{1})}+CT^{\frac12}\|\bar u_k\|_{C^{0}(0,T; C^{1})}\\
&\leq CT^{\frac12}\|\bar u_k\|_{C^{0}(0,T; C^{2})}+CT^{\frac12}\|\bar v_k\|_{C^{0}(0,T; C^{1})}+C\|\bar n_k\|_{L^{2}(0,T; L^{2})}.
\end{split}
\eeq 
Since $\bar n_k(\tau,0)=0$, we obtain that 
\beq\notag
\|\bar n_k\|_{L^{2}(0,T; L^{2})}\leq CT\|\bar n_k\|_{W^{2,1}_2([0,1]\times (0,T))}.
\eeq
If we choose $T>0$ small enough, we obtain
\beq\label{Lipn}
\begin{split}
\|\bar n_k\|_{W^{2,1}_2([0,1]\times (0,T))}
\leq & C(M_1, M_2, T)T^{\frac12}\left(\|\bar u_k\|_{C^{0}(0,T; C^{2})}+\|\bar v_k\|_{C^{0}(0,T; C^{1})}\right).
\end{split}
\eeq

\bigskip
\noindent \textbf{Step 4}. To obtain the estimates for $u_k$ and $v_k$, first notice that the equation of $u_k$ and $v_k$ can be understood in the weak senses, i.e., for any $\phi(x)\in \mathcal X_k$ and $t\in [0,T]$, it holds
\beq\label{applceL3}
\int_0^1\rho_k u_k\phi-\int\rho_0 u_0^k\phi=\int_0^t\int_0^1 \mathcal P^1(\rho_k, u_k, v_k, n_k)\phi+(\alpha_2+\alpha_3)\int_0^t\int_0^1 \dot n\cos n\sin n\phi_x,
\eeq
\beq\label{applceL4}
\int_0^1\rho_k v_k\phi-\int\rho_0 v_0^k\phi=\int_0^t\int_0^1 \mathcal P^2(\rho_k, u_k, v_k, n_k)\phi-\int_0^t\int_0^1 \big(\alpha_2\dot n_k\cos^2 n_k-\alpha_3\dot n_k\sin^2 n_k\big)\phi_x,
\eeq
where 
\beq\notag
\begin{split}
&\mathcal P^1(\rho_k, u_k, v_k, n_k)\\
=&(\alpha_0+\alpha_5+\alpha_6+\alpha_8)\big((u_k)_x\cos^2 n_k\big)_x+\alpha_1\big((u_k)_x\cos^4n_k\big)_x+(\alpha_4+\alpha_7)(u_k)_{xx}\\
&+\alpha_0\big((v_k)_x\cos n_k\sin n_k\big)_x+\alpha_1\big((v_k)_x\cos^3n_k\sin n_k\big)_x\\
&+\frac12(\alpha_2+\alpha_3+\alpha_5+\alpha_6)\big((v_k)_x\cos n_k\sin n_k\big)_x-(\rho_k u_k u_k)_x-\big(\rho_k^{\gamma}\big)_x-(n_k)_{xx}(n_k)_x,
\end{split}
\eeq
and 
\beq\notag
\begin{split}
\mathcal P^2(\rho_k, u_k, v_k, n_k)&=\alpha_1\big((v_k)_x\cos^2n_k\sin^2 n_k\big)_x+\frac12(-\alpha_2
+\alpha_5)\big((v_k)_x\cos^2 n_k\big)_x\\
&+\frac12(\alpha_3+\alpha_6)\big((v_k)_x\sin^2 n_k\big)_x+\frac12\alpha_4(v_k)_{xx}\\
&+\alpha_1\big((u_k)_x\cos^3n_k\sin n_k\big)_x+(\alpha_6+\alpha_8)\big((u_k)_x\cos n_k\sin n_k\big)_x-(\rho_k v_k v_k)_x.
\end{split}
\eeq
Similarly to the energy inequality \eqref{enest1}, we can obtain the
same form of energy estimates for the system \eqref{applce} so that
\beq\label{M3}
\|u_k\|_{C^0(0,T;C^2)}+\|v_k\|_{C^0(0,T;C^2)}
\leq C\|u_k\|_{C^0(0,T;L^2)}+C\|v_k\|_{C^0(0,T;L^2)}\leq M_3,
\eeq
provided $\inf\limits_{(x,t)}\rho_k(x,t)>0$. Here we have used the fact that the dimension of $\mathcal X_k$ is finite.

To apply the contraction map theorem, 
we define the linear operator $\mathcal N[\rho_k]:\, \mathcal X_k\rightarrow \mathcal X_k^*$ by
\beq\notag
\langle\mathcal N[\rho_k]\psi, \, \phi\rangle=\int_0^1\rho_k\psi\phi, \ \psi, \phi\in \mathcal X_k.
\eeq
It is easy to see that 
\beq\notag
\|\mathcal N[\rho_k]\|_{\mathcal L(\mathcal X_k, \mathcal X_k^*)}
\leq C(k)\|\rho_k\|_{L^1}.
\eeq
If $\inf\limits_{x}\rho_k>0$, the operator $\mathcal N[\rho_k]$ is invertible and 
\beq\notag
\|\mathcal N^{-1}[\rho_k]\|_{\mathcal L(\mathcal X_k^*, \mathcal X_k)}
\leq \left(\inf\limits_{x}\rho_k\right)^{-1}.
\eeq
Furthermore, for any $\rho_k^i\in L^1$ and $\inf\limits_{x}\rho_k^i>0$, $i=1,2$, it is easy to see that 
\beq\notag
\mathcal N^{-1}[\rho_k^1]-\mathcal N^{-1}[\rho_k^2]
=\mathcal N^{-1}[\rho_k^2]\left(\mathcal N[\rho_k^2]-\mathcal N[\rho_k^1]\right)\mathcal N^{-1}[\rho_k^1]
\eeq
which implies that
\beq\label{LipN}
\left\|\mathcal N^{-1}[\rho_k^1]-\mathcal N^{-1}[\rho_k^2]\right\|_{\mathcal L(\mathcal X_k^*, \mathcal X_k)}
\leq C\left\|\mathcal N[\rho_k^1]-\mathcal N[\rho_k^2]\right\|_{\mathcal L(\mathcal X_k, \mathcal X_k^*)}
\leq C\|\rho_k^1-\rho_k^2\|_{L^1}.
\eeq

Hence by the estimates \eqref{Liprho}, \eqref{Lipn} and \eqref{LipN}, we can apply the standard contraction map theorem to obtain the local existence of a unique solution $u_k, v_k\in C(0,T_k; \mathcal X_k)$ to \eqref{applceL3} and \eqref{applceL4} for some $T_k>0$. Then by the equations \eqref{applceL1} and \eqref{applceL2}, we can solve for $\rho_k, n_k$, which provides a unique local solution to the approximated system \eqref{applce} for
any fixed $k$.

\bigskip
\noindent \textbf{Step 5}. In this step, we will establish a uniform estimate of the local solution until $T_k$ in order to extend the solution beyond $T_k$ to any time $T>0$, which implies the existence of unique global solution of the system \eqref{applce} for any fixed $k$. We first show the following uniform estimate for $\rho_k$

\medskip
\noindent{\it Claim: For any $x\in [0,1]$ and $t\in [0,T_k]$, it holds
\beq\label{ufmrhok}
\frac{1}{c_1e^{t}}\leq\rho_k(x,t)\leq c_1e^{t}
\eeq
for some constant $c_1>0$.}

\medskip
\noindent Indeed, similar to the energy inequality \eqref{enest1}, we can obtain the same form of energy estimate for system \eqref{applce} so that 
\beq\label{M4}
\begin{split}
&\|(u_k)_x\|_{L^2(0,T_k;H^2)}+\|(v_k)_x\|_{L^2(0,T_k;H^2)}\\
\leq& C\|(u_k)_x\|_{L^2((0,1)\times(0,T_k)})+C\|(v_k)_x\|_{L^2((0,1)\times(0,T_k))}\leq M_4.
\end{split}
\eeq
By the first equation of \eqref{applce}, we can find $x_0(t)\in (0,1)$ such that 
\beq\notag
\rho_k(x_0(t),t)=\int_0^1\rho_k=\int_0^1\rho_0=1.
\eeq
Then
\beq\notag
\begin{split}
\frac{1}{\rho_k(x,t)}
=&\frac{1}{\rho_k(x_0(t),t)}+\int_{x_0(t)}^x\left(\frac{1}{\rho_k}\right)_y
\leq 1+\frac12\left\|\frac{1}{\rho_k(x,t)}\right\|_{L^\infty}
+\frac12\int_0^1\rho_k\left|\left(\frac{1}{\rho_k}\right)_x\right|^2
\end{split}
\eeq
which implies
\beq\label{pfclm1}
\begin{split}
\left\|\frac{1}{\rho_k(x,t)}\right\|_{L^\infty}
\leq 2+\int_0^1\rho_k\left|\left(\frac{1}{\rho_k}\right)_x\right|^2.
\end{split}
\eeq
By the first equation of \eqref{applce}, we have 
\beq\label{pfclm3}
\begin{split}
&\frac{d}{dt}\int_0^1\rho_k\left|\left(\frac{1}{\rho_k}\right)_x\right|^2\\
=&\int_0^1(\rho_k)_t\left|\left(\frac{1}{\rho_k}\right)_x\right|^2+2\int_0^1\rho_k\left(\frac{1}{\rho_k}\right)_x\left(\frac{1}{\rho_k}\right)_{xt}\\
=&-\int_0^1(\rho_ku_k)_x\left|\left(\frac{1}{\rho_k}\right)_x\right|^2+2\int_0^1\rho_k\left(\frac{1}{\rho_k}\right)_x\left(\frac{(\rho_ku_k)_x}{\rho_k^2}\right)_{x}.
\end{split}
\eeq
The last term on the right hand side can be computed by 
\beq\label{pfclm4}
\begin{split}
&2\int_0^1\rho_k\left(\frac{1}{\rho_k}\right)_x\left(\frac{(\rho_ku_k)_x}{\rho_k^2}\right)_{x}\\
=&2\int_0^1\rho_k\left(\frac{1}{\rho_k}\right)_x\left[\left(\left(-\frac{1}{\rho_k}\right)_xu_k\right)_x+\left(\frac{(u_k)_x}{\rho_k}\right)_x\right]\\
=&-\int_0^1\rho_ku_k\frac{\partial}{\partial x}\left|\left(\frac{1}{\rho_k}\right)_x\right|^2+2\int_0^1\left(\frac{1}{\rho_k}\right)_x(u_k)_{xx}.
\end{split}
\eeq
Combining \eqref{pfclm4} with \eqref{pfclm3}, we conclude that 
\beq\label{pfclm5}
\frac{d}{dt}\int_0^1\rho_k\left|\left(\frac{1}{\rho_k}\right)_x\right|^2=2\int_0^1\left(\frac{1}{\rho_k}\right)_x(u_k)_{xx}.
\eeq
The right hand side can be estimated as follows
\beq\notag
\begin{split}
&\left|\int_0^1\left(\frac{1}{\rho_k}\right)_x(u_k)_{xx}\right|\\
\leq&\int_0^1\rho_k^{\frac12}\left|\left(\frac{1}{\rho_k}\right)_x\right|\rho_k^{-\frac12}|(u_k)_{xx}|\\
\leq &\frac12\int_0^1\rho_k\left|\left(\frac{1}{\rho_k}\right)_x\right|^2+\frac12\left\|\frac{1}{\rho_k(x,t)}\right\|_{L^\infty}\int_0^1|(u_k)_{xx}|^2\\
\leq&\frac12\left(1+\int_0^1|(u_k)_{xx}|^2\right)\int_0^1\rho_k\left|\left(\frac{1}{\rho_k}\right)_x\right|^2+\int_0^1|(u_k)_{xx}|^2.
\end{split}
\eeq
where we have used \eqref{pfclm1} in last inequality. Denote 
$$
\mathcal Q(\rho_k)=\int_0^1\rho_k\left|\left(\frac{1}{\rho_k}\right)_x\right|^2, \quad a(t)=1+\int_0^1|(u_k)_{xx}|^2.
$$
Then by \eqref{pfclm4},  we have 
\beq\notag\frac{d}{dt}\mathcal Q(\rho_k)\leq a(t)\mathcal Q(\rho_k)
+\int_0^1|(u_k)_{xx}|^2
\eeq
which is equivalent to 
\beq\notag
\mathcal Q(\rho_k)-\mathcal Q(\rho_0)\leq 
2\int_0^t\int_0^1|(u_k)_{xx}|^2+\int_0^t a(t)\mathcal Q(\rho_k)
\leq 2M_4+\int_0^t a(t)\mathcal Q(\rho_k),
\eeq
where we have used \eqref{M4} in last step. 
By the Gronwall inequality, we obtain 
\beq\label{pfclm2}
\mathcal Q(\rho_k)
\leq \big(\mathcal Q(\rho_0)+2M_4\big)\exp\left(\int_0^ta(s)\right)
\leq \big(\mathcal Q(\rho_0)+2M_4\big)\exp\left(t+M_4\right)\leq Ce^t.
\eeq
Combining \eqref{pfclm1} and \eqref{pfclm2} together, we can prove the left hand side of \eqref{ufmrhok}. 

Denote $\gamma=1+2\delta$ for some $\delta>0$. Then it holds
\beq\notag
\begin{split}
\|\rho_k^\delta\|_{L^\infty}
\leq \int_0^1 \rho_k^{\delta}+\delta\int_0^1\rho_k^{\delta-1}(\rho_k)_x
\leq \left(\int_0^1 \rho_k^{\gamma}\right)^{\frac{\delta}{\gamma}}+C\left(\int_0^1 \rho_k^{\gamma}\right)^{\frac12}\left(\mathcal Q(\rho_k)\right)^{\frac12}\leq Ce^t,
\end{split}
\eeq
which completes the proof of the Claim.

\medskip
By using the uniform estimate \eqref{ufmrhok} for $\rho_k$ and the energy inequality, we can show 
\beq\notag
\|u_k\|_{C^0(0,T_k;\mathcal X_k)}+\|v_k\|_{C^0(0,T_k;\mathcal X_k)}
\leq C\|u_k\|_{C^0(0,T_k;L^2)}+C\|v_k\|_{C^0(0,T_k;L^2)}\leq M_5.
\eeq
Therefore, we can extend the solution beyond $T_k$ to any time $T>0$, which implies the existence of 
a unique smooth solution of the system \eqref{applce} for any fixed $k$.


\section{Existence of global weak solutions}

\noindent \textbf{Step 1.} 
Taking $k\rightarrow \infty$ in the approximated system \eqref{applce}, we may obtain the existence of 
a global weak solution with a smooth initial and boundary value and $\rho_0>\delta>0$. Since the limit process of this step is similar to the next step when $\delta\rightarrow 0$, we omit the details of this step.

\bigskip
\noindent \textbf{Step 2.} We first approximate the general initial and boundary data in Theorem \ref{mainth1} by smooth functions. 
We may extend $n$ to $\tilde n_0\in H^1(\R)$ such that $n_0=\tilde n_0$ on $(0,1)$, and obtain the smooth approximation of initial data by the standard mollification as follows
\beq\notag
\rho_0^\delta=\eta_\delta*\hat\rho_0+\delta,\quad 
u_0^\delta=\frac{1}{\sqrt{\rho_0^\delta}}\eta_{\delta}*\left(\widehat{\frac{m_0}{\sqrt{\rho_0}}}\right),\quad
v_0^\delta=\frac{1}{\sqrt{\rho_0^\delta}}\eta_{\delta}*\left(\widehat{\frac{l_0}{\sqrt{\rho_0}}}\right),\quad
n_0^\delta=\frac{\eta_\delta*\tilde n_0}{\big|\eta_\delta*\tilde n_0\big|}
\eeq
where, for small $\delta>0$, $\eta_\delta=\frac{1}{\delta}\eta\left(\frac{\cdot}{\delta}\right)$ is the standard mollifier, $\hat f$ is the zero extension of $f$ from $(0,1)$ to $\R$. Therefore $\rho_0^\delta, u_0^\delta, v_0^\delta, n_0^\delta \in C^{2+\alpha}([0,1])$ for $0<\alpha<1$, and it holds
\beq\label{appinitial1}
\rho_0^\delta\geq	\delta>0, \quad \rho_0^\delta\rightarrow \rho_0 \mbox{ in } L^{\gamma}, \quad  n_0^\delta\rightarrow n_0 \mbox{ in } H^{1},
\eeq
\beq\label{appinitial2}
 \sqrt{\rho_0^\delta}u_0^\delta\rightarrow \frac{m_0}{\sqrt{\rho_0}} \mbox{ in } L^{2}, \quad  \sqrt{\rho_0^\delta}v_0^\delta\rightarrow \frac{l_0}{\sqrt{\rho_0}} \mbox{ in } L^{2},\quad
 \rho_0^\delta u_0^\delta\rightarrow m_0 \mbox{ in } L^{\frac{2\gamma}{\gamma+1}},\quad \rho_0^\delta v_0^\delta\rightarrow l_0 \mbox{ in } L^{\frac{2\gamma}{\gamma+1}}
\eeq
as $\delta\rightarrow 0$.

Let $(\rho_\delta, u_\delta, v_\delta, n_\delta)$ be a sequence of global weak solutions to
\begin{equation}\label{lcedelta}
\begin{cases}
(\rho_\delta)_t+(\rho_\delta u_\delta)_x=0,\quad \rho_\delta>0,\\ 
(\rho_\delta u\delta)_t+(\rho_\delta u_\delta ^2)_x+\big(\rho_\delta^{\gamma}\big)_x=J^1_\delta-(n_\delta)_{xx}(n_\delta)_x,
\\
(\rho_\delta v_\delta)_t+(\rho_\delta u_\delta v_\delta)_x=J^2_\delta,\\
\gamma_1\left(\dot n_\delta-\frac12(v_\delta)_x\right)
-\gamma_2\left((u_\delta)_x\cos n_\delta\sin n_\delta+\frac12(v_\delta)_x(1-2\cos^2 n_\delta)\right) =(n_\delta)_{xx},
\end{cases}
\end{equation}
with the initial and boundary values 
\beq\label{deltainitial}
(\rho_\delta,\,  u_\delta,\, v_\delta,\, n_\delta)(x,0)=(\rho_0^\delta,\, u_0^\delta,\, v_0^\delta,\, n_0^\delta)(x),
\eeq
\beq\label{deltabdyvalue}
u_\delta(0,t)=v_\delta(0,t)=u_\delta(1,t)=v_\delta(1,t)=0,\quad (n_\delta)_x(0,t)=(n_\delta)_x(1,t)=0.
\eeq
Here $J^1_\delta$ and $J^2_\delta$ have the same forms as $J^1$ and $J^2$, but with $(u,v,n)$ replaced by $(u_\delta, v_\delta, n_\delta)$.

By Lemma \ref{lemma1}--Lemma  \ref{lemma3}, we can find a subsequence $(\rho_{\delta}, u_{\delta}, v_{\delta}, n_{\delta})$, still denoted as $(\rho_{\delta}, u_{\delta}, v_{\delta}, n_{\delta})$, such that for any $T>0$,
as $\delta\rightarrow 0$,
\beq\label{pfsec5.1}
\rho_{\delta}\overset{*}{\rightharpoonup} \rho, \mbox{ in } L^\infty(0,T;L^{\gamma}), \qquad 
\rho_{\delta}\rightharpoonup \rho,\mbox{ in } L^{2\gamma}([0,1]\times[0,T]),
\eeq
\beq\label{pfsec5.2}
\rho_{\delta}^{\gamma}\rightharpoonup \overline{\rho^{\gamma}},\mbox{ in } L^{2}([0,1]\times[0,T]),
\eeq
\beq\label{pfsec5.3}
u_{\delta}\rightharpoonup u,\mbox{ in } L^{2}(0,T; H_0^1),\quad
v_{\delta}\rightharpoonup v,\mbox{ in } L^{2}(0,T; H_0^1),
\eeq
\beq\label{pfsec5.4}
n_{\delta}\overset{*}{\rightharpoonup} n,\mbox{ in } L^{\infty}([0,1]\times[0,T]),\quad
(n_{\delta})_x\overset{*}{\rightharpoonup} n_x,\mbox{ in } L^{\infty}(0,T; L^2),
\eeq
\beq\label{pfsec5.5}
(n_{\delta})_t\rightharpoonup n_t,\mbox{ in } L^{2}([0,1]\times[0,T]),
\quad 
(n_{\delta})_{xx}\rightharpoonup n_{xx},\mbox{ in } L^{2}([0,1]\times[0,T]).
\eeq
Since $\rho_{\delta}>0$, for any nonnegative function $f\in C^\infty_0((0,1)\times (0,T))$ 
it holds that 
\beq\notag
\int_0^T\int_0^1 \rho f
=\lim\limits_{\delta\rightarrow 0}\int_0^T\int_0^1 \rho_{\delta} f\geq 0.
\eeq 
Since $f$ is arbitrary, we conclude that $\rho\geq 0$ a.e. in $(0,1)\times (0,T)$.

We need to show the limit $(\rho, u, v, n)$ is a solution to the system \eqref{applce}. We first state several compactness results that will be used in our proof. 

\begin{lemma}(\cite{simon90})\label{lemma5.1}
Assume $X\subset E\subset Y$ are Banach spaces and $X \hookrightarrow\hookrightarrow E$
is compact. Then the following embeddings are compact
$$
\left\{f:\, f\in L^{q}(0,T;X),\, \frac{\partial f}{\partial t}\in L^1(0,T; Y)\right\}\hookrightarrow\hookrightarrow L^q(0,T;E), \mbox{ for any } 1\leq q\leq \infty,
$$ 
$$
\left\{f:\, f\in L^{\infty}(0,T;X),\, \frac{\partial f}{\partial t}\in L^r(0,T; Y)\right\}\hookrightarrow\hookrightarrow C([0,T];E), \mbox{ for any } 1<r<\infty.
$$ 
\end{lemma}

\begin{lemma}(\cite{feireisl04})\label{lemma5.2}
Let $\bar O\subset \mathbb R^n$ be compact and $X$ be a separable Banach space. Assume that $f_{\delta}:\bar O\rightarrow X^*$ is a sequence of measurable functions such that for any $k$
$$
\mbox{ess}\sup\limits_{\bar O}\|f_{\delta}\|_{X^*}\leq N<\infty.
$$
Moreover, the family of functions $\langle f_{\delta}, \Phi \rangle$ is equi-continuous for any $\Phi$ belonging to a dense subset of $X$. Then $f_{\delta}\in C(\bar O; X-w)$ for any $k$, i.e., for any $g\in X*$, $\langle f_{\delta},\, g\rangle\in C(\bar O)$. Furthermore, there exists $f\in  C(\bar O; X-w)$ such that (after taking possible subsequences)
$$
f_{\delta}\rightarrow f,\quad \mbox{in }  C(\bar O; X-w)
$$
as $\delta\rightarrow 0$.
\end{lemma}

First observe that $\rho_{\delta}\in L^{2\gamma}([0,1]\times[0,T])$ and $u_{\delta}\in L^2(0,T;H^1_0)\subset L^2(0,T;L^{\infty})$ imply
\beq\notag
\rho_{\delta}u_{\delta}\in L^{\frac{2\gamma}{\gamma+1}}(0,T;L^{2\gamma}),\quad 
(\rho_{\delta})_t=-(\rho_{\delta}u_{\delta})_x\in L^{\frac{2\gamma}{\gamma+1}}(0,T;H^{-1}).
\eeq 
By Lemma \ref{lemma5.1} and Lemma \ref{lemma5.2}, and $\frac{2\gamma}{\gamma+1}>1$, $\rho_{\delta}\in L^\infty(0,T;L^{\gamma})$, $L^{\gamma} \hookrightarrow\hookrightarrow H^{-1}$, we conclude 
\beq\label{pfsec5.6}
\rho_{\delta}\rightarrow \rho, \mbox{ in } C(0,T; L^{\gamma}-\omega),\quad 
\rho_{\delta}\rightarrow \rho, \mbox{ in } C(0,T; H^{-1}),
\eeq 
where $f\in C(0,T; X-\omega)$ if for any $g\in X^*$, $\langle f(t),\, g\rangle\in C([0,T])$. Hence 
\beq\label{pfsec5.12}
\rho_{\delta}u_{\delta}\rightarrow \rho u,\mbox{ in }\mathcal D'((0,1)\times(0,T)),\quad 
\rho_{\delta}v_{\delta}\rightarrow \rho v,\mbox{ in }\mathcal D'((0,1)\times(0,T)),
\eeq
and furthermore
\beq\label{pfsec5.7}
\rho_t+(\rho u)_x=0,\mbox{ in }\mathcal  D'((0,1)\times(0,T)).
\eeq
By \eqref{pfsec5.6}, it also holds that
\beq\label{pfsec5.8}
\rho(x,0)=\rho_0(x),\mbox{ weakly in } L^{\gamma}([0,1]).
\eeq

By the fact $(n_{\delta})_t \in L^2(0,T;L^2)$, \eqref{pfsec5.4} and \eqref{pfsec5.5}, we can apply Lemma \ref{lemma5.1} to obtain
 \beq\label{pfsec5.9}
n_{\delta}\rightarrow n, \mbox{ in } C([0,1]\times[0,T]),\quad 
n_{\delta}\rightarrow n, \mbox{ in } L^2(0,T; C^{1}),
\eeq 
Combining with \eqref{pfsec5.3}-\eqref{pfsec5.5}, we can show the limit $n$ satisfies the following equation:
\beq\label{pfsec5.10}
\gamma_1\left(\dot n-\frac12v_x\right)
-\gamma_2\left(u_x\cos n\sin n+\frac12v_x(1-2\cos^2 n)\right) =n_{xx}.
\eeq
By \eqref{pfsec5.9}, it also holds that
\beq\label{pfsec5.11}
n(x,0)=n_0(x),\mbox{ in } [0,1].
\eeq

By the fact $\sqrt{\rho_{\delta}}\in L^{2\gamma}([0,1]\times[0,T])$ and $\sqrt{\rho_{\delta}}u_{\delta}\in L^{\infty}(0,T;L^2)$, it holds
\beq\notag
\rho_{\delta}u_{\delta}\in L^{\infty}(0,T;L^{\frac{2\gamma}{\gamma+1}}).
\eeq 
Combining with \eqref{pfsec5.3}, we have
\beq\label{pfsec5.13} 
\rho_{\delta}u_\delta^2\rightharpoonup \rho u^2, \mbox{ in } L^2(0,T; L^{\frac{2\gamma}{\gamma+1}}).
\eeq
By the second equation of system \eqref{lcedelta}, we have 
\beq\notag
(\rho_\delta u_\delta)_t=-(\rho_\delta u_\delta^2)_x-\big(\rho_\delta^{\gamma}\big)_x+J^1_\delta-(n_\delta)_{xx}(n_\delta)_x \in L^2(0,T;W^{-1, \frac{2\gamma}{\gamma+1}}),
\eeq
where $J^1_\delta$ has the same form as $J^1$, but with $(u,v,n)$ replaced by $(u_\delta, v_\delta, n_\delta)$.
By using  Lemma \ref{lemma5.1} and Lemma \ref{lemma5.2}, we conclude 
\beq\label{pfsec5.14}
\rho_{\delta}u_\delta\rightarrow \rho u, \mbox{ in } C(0,T; L^{\frac{2\gamma}{\gamma+1}}-\omega),\quad 
\rho_{\delta}u_\delta\rightarrow \rho u, \mbox{ in } C(0,T; H^{-1}).
\eeq 
Combining with \eqref{pfsec5.3}, we conclude that 
\beq\label{pfsec5.15}
\rho_{\delta}u_{\delta}^2\rightarrow \rho u^2,\mbox{ in }\mathcal  D'((0,1)\times(0,T)).
\eeq
Therefore
\beq\label{pfsec5.16}
(\rho u)_t+(\rho u^2)_x+\big(\overline{\rho^{\gamma}}\big)_x=J^1-n_{xx}n_x,\mbox{ in }\mathcal  D'((0,1)\times(0,T)).
\eeq
By \eqref{pfsec5.14}, it holds that
\beq\label{pfsec5.17}
\rho u(x,0)=m_0(x),\mbox{ weakly in } L^{\frac{2\gamma}{\gamma+1}}([0,1]).
\eeq
Similarly, we can also prove that 
\beq\label{pfsec5.18}
(\rho v)_t+(\rho u v)_x=J^2,\mbox{ in }\mathcal  D'((0,1)\times(0,T)),
\eeq
\beq\label{pfsec5.19}
\rho v(x,0)=n_0(x),\mbox{ weakly in } L^{\frac{2\gamma}{\gamma+1}}([0,1]).
\eeq

By \eqref{pfsec5.15}, for some $t\in (0,T)$ and small $\epsilon>0$, it holds
\beq\notag
\frac{1}{\epsilon}\int_t^{t+\epsilon}\int_0^1\rho u^2
=\frac{1}{\epsilon}\int_t^{t+\epsilon}\lim\limits_{\delta\rightarrow 0}\int_0^1\rho_\delta u^2_\delta
\leq \frac{1}{\epsilon}\int_t^{t+\epsilon}\overline{\lim\limits_{\delta\rightarrow 0}}\int_0^1\rho_\delta u^2_\delta.
\eeq 
Sending $\epsilon\rightarrow 0^+$ and using the Lebesgue Differentiation Theorem, we obtain
\beq\notag
\int_0^1\rho u^2
\leq \overline{\lim\limits_{\delta\rightarrow 0}}\int_0^1\rho_\delta u^2_\delta,
\eeq
for a.e. $t\in (0,T)$. Combining this limit with the lower semicontinuity, we can prove the energy inequality is valid. 

\bigskip
The only thing left is to show
$\overline{\rho^\gamma}=\rho^\gamma$. To this end, we denote 
\beq\notag
A(n)=(A_{ij}(n))_{2\times 2}
\eeq
where the elements of $A_{ij}$ are given as follows
\beq\notag
A_{11}(n)=(\alpha_0+\alpha_5+\alpha_6+\alpha_8)\cos^2 n+\alpha_1\cos^4n+(\alpha_4+\alpha_7)
\eeq
\beq\notag
A_{12}(n)= \alpha_0\cos n\sin n+\alpha_1\cos^3n\sin n+\frac12(\alpha_2+\alpha_3+\alpha_5+\alpha_6)\cos n\sin n
\eeq
\beq\notag
A_{21}(n)=\alpha_1\cos^3n\sin n+(\alpha_6+\alpha_8)\cos n\sin n
\eeq
\beq\notag
A_{22}(n)=\alpha_1\cos^2n\sin^2 n+\frac12(-\alpha_2+\alpha_5)\cos^2 n+\frac12(\alpha_3+\alpha_6)\sin^2 n+\frac12\alpha_4.
\eeq
By the relations \eqref{alphas}, direct computations imply that there exist two positive constants $\lambda,\Lambda<\infty$ such that for any $\y\in \R^2$
\beq \label{postiveA}
\lambda|\y|^2\leq \y^TA(n)\y\leq \Lambda|\y|^2.
\eeq
In fact
\beq\notag
\begin{split}
\y^TA(n)\y=&A_{11}(n)y_1^2+(A_{12}(n)+A_{21}(n))y_1y_2+A_{22}(n)y_2^2\\
=&\big[(\alpha_0+\alpha_5+\alpha_6+\alpha_8)\cos^2 n+\alpha_1\cos^4n+(\alpha_4+\alpha_7)\big]y_1^2\\
&+\big[(\alpha_0+\alpha_6+\alpha_8)\cos n\sin n+2\alpha_1\cos^3n\sin n+\frac12(\alpha_2+\alpha_3+\alpha_5+\alpha_6)\cos n\sin n\big]y_1y_2\\
&+\big[\alpha_1\cos^2n\sin^2 n+\frac12(-\alpha_2+\alpha_5)\cos^2 n+\frac12(\alpha_3+\alpha_6)\sin^2 n+\frac12\alpha_4\big]y_2^2\\
=&\frac14\left(
 \frac{\gamma_2}{\sqrt{\gamma_1}}y_1\sin(2n)+ \frac{1}{\sqrt{\gamma_1}}(\gamma_1-\gamma_2\cos(2n))y_2
\right)^2\\
&+\frac14\left(-\alpha_1-\frac{\gamma_2^2}{\gamma_1}\right)y_1^2+(\alpha_4+\alpha_7)y_1^2
+\frac14\left(2\alpha_4+\alpha_5+\alpha_6-\frac{\gamma_2^2}{\gamma_1}\right)y_2^2\\
&\frac14\left(\alpha_1+\frac{\gamma_2^2}{\gamma_1}\right)\left(y_1\cos(2n)+y_2\sin(2n)\right)^2\\
&+(\alpha_0+\alpha_1+\alpha_5+\alpha_6+\alpha_8)\left[\big(y_1\cos n+\frac12y_2\sin n\big)^2-\frac14 y_2^2\sin^2 n\right].
\end{split}
\eeq
Therefore
\beq\notag
\begin{split}
\y^TA(n)\y\geq& \frac14\left(-\alpha_1-\frac{\gamma_2^2}{\gamma_1}\right)y_1^2+(\alpha_4+\alpha_7)y_1^2
+\frac14\left(2\alpha_4+\alpha_5+\alpha_6-\frac{\gamma_2^2}{\gamma_1}\right)y_2^2\\
&-\frac14(\alpha_0+\alpha_1+\alpha_5+\alpha_6+\alpha_8) y_2^2\sin^2 n.
\end{split}
\eeq
If we take 
\beq\notag
\lambda=\min\left\{
 (\alpha_4+\alpha_7)-\frac14\left(\alpha_1+\frac{\gamma_2^2}{\gamma_1}\right),\ 
 \left(2\alpha_4+\alpha_5+\alpha_6-\frac{\gamma_2^2}{\gamma_1}\right)-(\alpha_0+\alpha_1+\alpha_5+\alpha_6+\alpha_8)
\right\},
\eeq
then by the relation \eqref{alphas}, we know that $\lambda>0$ and we have shown the estimate \eqref{postiveA}.

By the definition of $A(n)$, we see that the matrix valued function $A(\cdot)\in C^{\infty}$. By the estimate \eqref{postiveA}, the inverse matrix function $A^{-1}$ exists and 
$$
\frac{d}{dn}\big(A^{-1}(n)\big)=A^{-1}\frac{d}{dn}\big(A(n)\big)A^{-1}.
$$

The equations for $\bu=(u,\, v)^T$ can be written as 
\beq\label{vec1dlce}
\rho\bu_t+\rho u\bu_x+ {\bf P}_x=\big(A(n)\bu_x\big)_x+(B_1(n))_x-B_2(n)
\eeq
where 
\beq\notag
{\bf P}=(\overline{\rho^{\gamma}},\, 0)^T,
\eeq
\beq\notag
B_1(n)=\big((\alpha_2+\alpha_3)\dot n\cos n\sin n,\ \alpha_2\dot n\cos^2 n-\alpha_3\dot n\sin^2 n\big)^T,
\eeq
\beq\notag
B_2(n)=\big(n_{xx}n_x,\, 0\big)^T.
\eeq
Similarly, we can 
rewrite the equations for $\bu_\delta=(u_\delta,\, v_\delta)^T$, ${\bf P}_\delta=(\rho_\delta^{\gamma},\, 0)^T$ in the similar form
\beq\label{Dvec1dlce}
\rho_\delta(\bu_\delta)_t+\rho_\delta u_\delta(\bu_\delta)_x+ ({\bf P_\delta})_x=\big(A(n_\delta)(\bu_\delta)_x\big)_x+(B_1(n_\delta))_x-B_2(n_\delta).
\eeq
Denote  
\beq\notag
\h=\bu_x-A^{-1}(n)\p,\quad \h_\delta=(\bu_\delta)_x-A^{-1}(n_\delta)\p_\delta.
\eeq
We have the following lemma. 

\begin{lemma}\label{lemma5.3}
As $\delta\rightarrow 0$, it holds
\beq\label{pfsec6.1}
\rho_{\delta}\h_\delta\rightarrow \rho \h,\mbox{ in }\mathcal  D'((0,1)\times(0,T)).
\eeq
\end{lemma}

\pf 
The main difficulty of the proof arises from $\rho u\not\in L^2$. To overcome it, we need to mollify the density $\rho$ by $\langle\hat\rho\rangle_\sigma=\eta_\sigma*\hat\rho$, where $\eta_\sigma=\frac{1}{\sigma}\sigma\left(\frac{\cdot}{\sigma}\right)$ is the standard mollifier, $\hat f$ is the zero extension of $f$ from $(0,1)$ to $\R$.
By Lemma 3.3 in \cite{feireisl04}, the zero-extension of $\hat\rho$ still satisfies the same equation
\beq\label{rhohat}
(\hat\rho)_t+(\hat\rho \hat u)_x=0, \quad \mbox{in }\mathcal D'(\R\times(0,T)).
\eeq
Denote $\tau^{\sigma}=(\rs\hat u)_x-\langle(\hat \rho\hat u)_x\rangle_\sigma$. By Lemma 2.3 in \cite{lions96}, we know that $\tau^{\sigma}\in L^{\frac{2\gamma}{\gamma+1}}(\R\times(0,T)$, and as $\sigma\rightarrow 0$
\beq\label{taus}
\tau^{\sigma}\rightarrow 0,\quad \mbox{in }L^1(\R\times(0,T)).
\eeq
Taking the standard mollifier as the test function, we obtain 
\beq\label{rhosigma}
(\rs)_t+(\rs \hat u)_x=\tau^\sigma, \quad \mbox{in }\mathcal D'(\R\times(0,T)).
\eeq
Similarly, it also hold for the approximate solutions
\beq\label{rhods}
(\rds)_t+(\rds \hat u_\delta)_x=\tau^\sigma_\delta, \quad \mbox{in }\mathcal D'(\R\times(0,T)),
\eeq
where $\tau^\sigma_\delta$ has the same form as $\tau^\sigma$, but with $\rho, u$ replaced by $\rho_\delta, u_\delta$.  We also know that, for any $\delta>0$, $\tau^{\sigma}_\delta\in L^{\frac{2\gamma}{\gamma+1}}(\R\times(0,T)$, and as $\sigma\rightarrow 0$
\beq\label{tauds}
\tau^{\sigma}_\delta\rightarrow 0,\quad \mbox{in }L^1(\R\times(0,T)).
\eeq
%

Multiplying the equation \eqref{Dvec1dlce} by $\varphi\phi A^{-1}(n_\delta)\int_0^x\rds$ from left for any $\varphi\in C^\infty_0(0,T)$ and $\phi\in C^\infty_0(0,1)$, and integrating by parts, we obtain
\beq\notag
\begin{split}
&\int_0^T\int_0^1\varphi\phi\h_\delta\rds\\
=&\int_0^T\int_0^1\varphi'\phi\rho_\delta  A^{-1}(n_\delta)\bu_\delta\int_0^x\rds+\int_0^T\int_0^1\varphi\phi\rho_\delta  A^{-1}(n_\delta)\bu_\delta\left(\int_0^x\rds\right)_t\\
&+\int_0^T\int_0^1\varphi\phi\rho_\delta  \big(A^{-1}(n_\delta)\big)_t\bu_\delta\int_0^x\rds+\int_0^T\int_0^1\varphi\phi'\rho_\delta u_\delta A^{-1}(n_\delta)\bu_\delta\int_0^x\rds\\
&+\int_0^T\int_0^1\varphi\phi\rho_\delta\rds u_\delta A^{-1}(n_\delta)\bu_\delta+\int_0^T\int_0^1\varphi\phi\rho_\delta u_\delta \big(A^{-1}(n_\delta)\big)_x\bu_\delta\int_0^x\rds\\
&+\int_0^T\int_0^1\varphi\phi A^{-1}(n_\delta)(B_1(n_\delta))_x\int_0^x\rds-\int_0^T\int_0^1\varphi\phi A^{-1}(n_\delta)B_2(n_\delta)\int_0^x\rds\\
&-\int_0^T\int_0^1\varphi\phi'\h_\delta\int_0^x\rds-\int_0^T\int_0^1\varphi\phi A^{-1}(n_\delta)\big(A(n_\delta)\big)_x\h_\delta\int_0^x\rds.
\end{split}
\eeq
The equation \eqref{rhods} implies
$$
\frac{\partial}{\partial t}\left(\int_0^x\rds\right)=-\rds \hat u_\delta+\tau^\sigma_\delta.
$$
Using this fact, we have
\beq\notag
\begin{split}
&\int_0^T\int_0^1\varphi\phi\h_\delta\rds\\
=&\int_0^T\int_0^1\varphi'\phi\rho_\delta  A^{-1}(n_\delta)\bu_\delta\int_0^x\rds+\int_0^T\int_0^1\varphi\phi\rho_\delta  A^{-1}(n_\delta)\bu_\delta\int_0^x\tau^\sigma_\delta\\
&+\int_0^T\int_0^1\varphi\phi\rho_\delta  \big(A^{-1}(n_\delta)\big)_t\bu_\delta\int_0^x\rds+\int_0^T\int_0^1\varphi\phi'\rho_\delta u_\delta A^{-1}(n_\delta)\bu_\delta\int_0^x\rds\\
&+\int_0^T\int_0^1\varphi\phi\rho_\delta u_\delta \big(A^{-1}(n_\delta)\big)_x\bu_\delta\int_0^x\rds\\
&+\int_0^T\int_0^1\varphi\phi A^{-1}(n_\delta)(B_1(n_\delta))_x\int_0^x\rds-\int_0^T\int_0^1\varphi\phi A^{-1}(n_\delta)B_2(n_\delta)\int_0^x\rds\\
&-\int_0^T\int_0^1\varphi\phi'\h_\delta\int_0^x\rds-\int_0^T\int_0^1\varphi\phi A^{-1}(n_\delta)\big(A(n_\delta)\big)_x\h_\delta\int_0^x\rds.
\end{split}
\eeq
By the Lebesgue Dominated Convergence theorem and \eqref{tauds}, we may take the limit $\sigma\rightarrow 0$ and get 
\beq\label{pfsec6.3}
\begin{split}
&\int_0^T\int_0^1\varphi\phi\h_\delta\rho_\delta\\
=&\int_0^T\int_0^1\varphi'\phi\rho_\delta  A^{-1}(n_\delta)\bu_\delta\int_0^x\rho_\delta
+\int_0^T\int_0^1\varphi\phi\rho_\delta  \big(A^{-1}(n_\delta)\big)_t\bu_\delta\int_0^x\rho_\delta\\
&+\int_0^T\int_0^1\varphi\phi'\rho_\delta u_\delta A^{-1}(n_\delta)\bu_\delta\int_0^x\rho_\delta+\int_0^T\int_0^1\varphi\phi\rho_\delta u_\delta \big(A^{-1}(n_\delta)\big)_x\bu_\delta\int_0^x\rho_\delta\\
&+\int_0^T\int_0^1\varphi\phi A^{-1}(n_\delta)(B_1(n_\delta))_x\int_0^x\rho_\delta-\int_0^T\int_0^1\varphi\phi A^{-1}(n_\delta)B_2(n_\delta)\int_0^x\rho_\delta\\
&-\int_0^T\int_0^1\varphi\phi'\h_\delta\int_0^x\rho_\delta-\int_0^T\int_0^1\varphi\phi A^{-1}(n_\delta)\big(A(n_\delta)\big)_x\h_\delta\int_0^x\rho_\delta.
\end{split}
\eeq
By the definition of $B_2(n_\delta)$ and integration by parts, we obtain
\beq\label{pfsec6.6}
\begin{split}
&-\int_0^T\int_0^1\varphi\phi A^{-1}(n_\delta)B_2(n_\delta)\int_0^x\rho_\delta\\
=&\frac12\int_0^T\int_0^1\varphi\phi' A^{-1}(n_\delta)\big(|(n_\delta)_x|^2,\, 0\big)^T\int_0^x\rho_\delta
+\frac12\int_0^T\int_0^1\varphi\phi\big(A^{-1}(n_\delta)\big)_x\big(|(n_\delta)_x|^2,\, 0\big)^T\int_0^x\rho_\delta\\
&+\frac12\int_0^T\int_0^1\varphi\phi \rho_\delta A^{-1}(n_\delta)\big(|(n_\delta)_x|^2,\, 0\big)^T.
\end{split}
\eeq
By the definition of $B_1(n_\delta)$, we obtain
\beq\notag
\begin{split}
&\int_0^T\int_0^1\varphi\phi A^{-1}(n_\delta)\big(B_1(n_\delta)\big)_x\int_0^x\rho_\delta\\
=&\int_0^T\int_0^1\varphi\phi \left(A^{-1}(n_\delta)B_1(n_\delta)\right)_x\int_0^x\rho_\delta
-\int_0^T\int_0^1\varphi\phi \left(A^{-1}(n_\delta)\right)_xB_1(n_\delta)\int_0^x\rho_\delta.
\end{split}
\eeq
It is not hard to see that there is a vector function $\mathcal F(n_\delta)$ (smooth in $n_\delta$) such that
\beq\notag
A^{-1}(n_\delta)B_1(n_\delta)=\mathcal F_t(n_\delta)+u_\delta\mathcal F_x(n_\delta).
\eeq
Then
\beq\label{pfsec6.7}
\begin{split}
&\int_0^T\int_0^1\varphi\phi A^{-1}(n_\delta)\big(B_1(n_\delta)\big)_x\int_0^x\rho_\delta\\
=&-\int_0^T\int_0^1\varphi'\phi \mathcal F_x(n_\delta)\int_0^x\rho_\delta-\int_0^T\int_0^1\varphi\phi' u_\delta\mathcal F_x(n_\delta)\int_0^x\rho_\delta\\
&-\int_0^T\int_0^1\varphi\phi \left(A^{-1}(n_\delta)\right)_xB_1(n_\delta)\int_0^x\rho_\delta.
\end{split}
\eeq
To estimate the second term on right side of \eqref{pfsec6.3}, we use $\varphi\phi n $ as the test function for the first equation of \eqref{lcedelta} to obtain
\beq\notag
\int_0^T\int_0^1\varphi\phi\rho_\delta (n_\delta)_t
=-\int_0^T\int_0^1\varphi'\phi\rho_\delta n_\delta
-\int_0^T\int_0^1\varphi\rho_\delta u_\delta (n_\delta\phi)_x.
\eeq
Similarly, it holds
\beq\notag
\int_0^T\int_0^1\varphi\phi\rho n_t
=-\int_0^T\int_0^1\varphi'\phi\rho n
-\int_0^T\int_0^1\varphi\rho u (n\phi)_x.
\eeq
Taking the difference, and using \eqref{pfsec5.1}, \eqref{pfsec5.12} and \eqref{pfsec5.9}, we have
\beq\label{pfsec6.8}
\rho_{\delta}(n_{\delta})_t\rightarrow \rho n_t,\mbox{ in }\mathcal  D'((0,1)\times(0,T)).
\eeq
Furthermore, since
$$
\int_0^x\rho_\delta\in L^\infty(0,T;W^{1,\gamma}), \quad
\frac{\partial}{\partial t}\left(\int_0^x\rho_\delta\right)=-\rho_\delta u_\delta\in L^\infty\left(0,T;L^{\frac{2\gamma}{\gamma+1}}\right)
$$
we obtain by Lemma \ref{lemma5.1} and \eqref{pfsec5.1}
\beq\label{pfsec6.4}
\int_0^x\rho_\delta\rightarrow \int_0^x\rho, \quad \mbox{in }C([0,1]\times[0,T]),\quad \mbox{as }\delta\rightarrow 0.
\eeq

Now, we are ready to take limit in \eqref{pfsec6.3}.
Letting $\delta\rightarrow 0$ in \eqref{pfsec6.3} \eqref{pfsec6.6} and \eqref{pfsec6.7}, and using the facts \eqref{pfsec6.4}, \eqref{pfsec6.8}, \eqref{pfsec5.1}-\eqref{pfsec5.3}, \eqref{pfsec5.12}, \eqref{pfsec5.9} and \eqref{pfsec5.15}, we obtain
\beq\label{pfsec6.5}
\begin{split}
&\lim\limits_{\delta\rightarrow 0}\int_0^T\int_0^1\varphi\phi\h_\delta\rho_\delta\\
=&\int_0^T\int_0^1\varphi'\phi\rho  A^{-1}(n)\bu\int_0^x\rho
+\int_0^T\int_0^1\varphi\phi\rho  \big(A^{-1}(n)\big)_t\bu\int_0^x\rho\\
&+\int_0^T\int_0^1\varphi\phi'\rho u A^{-1}(n)\bu\int_0^x\rho+\int_0^T\int_0^1\varphi\phi\rho u \big(A^{-1}(n)\big)_x\bu\int_0^x\rho\\
&+\int_0^T\int_0^1\varphi\phi A^{-1}(n)(B_1(n))_x\int_0^x\rho-\int_0^T\int_0^1\varphi\phi A^{-1}(n)B_2(n)\int_0^x\rho\\
&-\int_0^T\int_0^1\varphi\phi'\h\int_0^x\rho-\int_0^T\int_0^1\varphi\phi A(n)\big(A^{-1}(n)\big)_x\h\int_0^x\rho.
\end{split}
\eeq

We may go through the same arguments for $\rho$ and $u$, and show that right side of \eqref{pfsec6.5} is exactly
$$
\int_0^T\int_0^1\varphi\phi\h\rho,
$$
which completes the proof of the lemma.
\endpf

We also need the following result. 

\begin{lemma}(\cite{feireisl04})\label{lemma5.5}
Let $\bar O\subset \mathbb R^n$ be a measurable set and $f_k \in L^1(O;\R^N)$ for $k\in \mathbb Z_+$ such that 
$$
f_k\rightharpoonup f, \quad \mbox{in }\ L^1(O;\R^N).
$$
Let $\Phi:\R^N\rightarrow (-\infty,\infty]$ be a lower semi-continuous convex function such that 
$\Phi(f_k)\in L^1(O)$ for any $k$ and 
$$
\Phi(f_k)\rightharpoonup \overline{\Phi(f)}, \quad \mbox{in }\ L^1(O).
$$
Then
$$
\Phi(f)\leq \overline{\Phi(f)}, \quad a.e.\ \mbox{in }\ O.
$$ 
Moreover, if $\Phi$ is strictly convex on an open convex set $U\subset \R^N$ and 
$$
\Phi(f)= \overline{\Phi(f)}, \quad a.e.\ \mbox{in }\ O,
$$ 
then 
$$
f_k\rightarrow f, \quad \mbox{for }\ a.e.\ y\in \{y\in O\,|\, f(y)\in U\}.
$$
\end{lemma}

\medskip
The proof of Theorem \ref{mainth1} will be completed by the following Lemma.

\begin{lemma}\label{lemma5.4}
As $\delta\rightarrow 0$, it holds
\beq\label{pfsec6.2}
\lim\limits_{\delta\rightarrow 0}\int_0^T\int_0^1\rho_\delta\log(\rho_\delta)
=\int_0^T\int_0^1\rho \log\rho.
\eeq
\end{lemma}

\pf By Proposition 4.2 in \cite{fnp01}, if $\rho\in L^2((0,1)\times(0,T))$, $u\in L^2(0,T;H_0^1)$ solves the equation 
\beq\notag
\rho_t+(\rho u)_x=0, \quad \mbox{in } \mathcal D'((0,1)\times(0,T))
\eeq
then 
\beq\label{pssec8.1}
(b(\rho))_t+(b(\rho)u)_x+(b'(\rho)\rho-b(\rho))u_x=0, \quad \mbox{in } \mathcal D'((0,1)\times(0,T))
\eeq
for any $b\in C^1(\mathbb R)$ such that $b'(x)\equiv 0$ for all large enough $x\in \mathbb R$.

For any positive integers $j, K$, we may take a family of functions $b_K^j\in C^1(\R)$ with
\beq\notag
b_K^j(x)=\left\{
\begin{array}{ll}
\displaystyle\left(x+\frac{1}{j}\right)\log\left(x+\frac{1}{j}\right),&\quad\mbox{if }0\leq x\leq K,\\
\displaystyle\left(K+1+\frac{1}{j}\right)\log\left(K+1+\frac{1}{j}\right),&\quad\mbox{if }x\geq K+1.
\end{array}
\right.
\eeq 
Since $\rho\in L^\infty(0,T;L^{\gamma})$, we have $\rho<\infty$ a.e. in $(0,1)\times(0,T)$. This implies that 
$b_K^j(\rho)\rightarrow (\rho+\frac1j)\log(\rho+\frac1j)$ a.e. in $(0,1)\times(0,T)$ as $K\rightarrow \infty$. Hence, by using the Lebesgue Dominated Convergence theorem, we conclude 
\beq\label{pfsec8.2}
\left(\left(\rho+\frac1j\right)\log\left(\rho+\frac1j\right)\right)_t+\left(\left(\rho+\frac1j\right)\log\left(\rho+\frac1j\right)u\right)_x+\Big(\rho -\frac1j\log\big(\rho+\frac1j\big)\Big)u_x=0, 
\eeq
in $\mathcal D'((0,1)\times(0,T))$. 

It is easy to see that $\left(\rho+\frac1j\right)\log\left(\rho+\frac1j\right)\in L^2((0,1)\times (0,T))$ since $\rho\in L^{2\gamma}((0,1)\times (0,T))$. By Lemma 3.3 in \cite{feireisl04}, the zero-extension of $\rho$ outside $(0,1)$ satisfies the same equation. By the mollification, the integration by parts and the limiting process, 
we may take the test function to be the constant $1$ so that
\beq\label{pfsec8.3}
\begin{split}
&\int_0^T\int_0^1\rho u_x\\
=&\int_0^1\left(\rho_0+\frac1j\right)\log\left(\rho_0+\frac1j\right)-\int_0^1\left(\rho+\frac1j\right)\log\left(\rho+\frac1j\right)(T)\\
&+\frac1j\int_0^T\int_0^1u_x\log\left(\rho+\frac1j\right).
\end{split}
\eeq

Similar estimates are valid for approximated solutions $\rho_\delta$, $u_\delta$. More precisely, we have
\beq\label{pfsec8.4}
\left(\rho_\delta\log\left(\rho_\delta\right)\right)_t+\left(\rho_\delta\log\left(\rho_\delta\right)u_\delta\right)_x+\rho_\delta (u_\delta)_x=0, 
\eeq
in $\mathcal D'((0,1)\times(0,T))$,
and  
\beq\label{pfsec8.5}
\begin{split}
\int_0^T\int_0^1\rho_\delta (u_\delta)_x
=\int_0^1\rho_0^\delta\log\left(\rho_0^\delta\right)
-\int_0^1\rho_\delta\log\left(\rho_\delta\right)(T)\\
\end{split}
\eeq
Since $\rho_\delta\in L^\infty(0,T;L^\gamma),$ we have
\beq\notag
\rho^\delta\log\left(\rho^\delta\right)\in L^\infty(0,T;L^{\tilde\gamma})
\eeq
for $1<\tilde{\gamma}<\gamma$. By the equation \eqref{pfsec8.4}, we obtain
\beq\notag
\left(\rho_\delta\log\left(\rho_\delta\right)\right)_t\in L^{\frac{2\gamma}{\gamma+1}}(0,T;W^{-1,\frac{2\gamma}{\gamma+1}}).
\eeq
By Lemma \ref{lemma5.2}, we conclude as $\delta\rightarrow 0$
\beq\notag
\rho^\delta\log\left(\rho^\delta\right)
\rightarrow 
\overline{\rho\log\left(\rho\right)},\quad \mbox{in } C([0,T]; L^{\tilde\gamma}-\omega).
\eeq
This implies 
\beq\label{pfsec8.10}
\lim\limits_{\delta\rightarrow 0}\int_0^1\rho^\delta\log\left(\rho^\delta\right)(T)
{\color{yellow}{=}}
\int_0^1\overline{\rho\log\left(\rho\right)}(T).
\eeq
Since the function $x\log\left(x\right)$ is convex for any $x>0$, 
Lemma \ref{lemma5.5} implies that 
\beq\label{pfsec8.6}
\rho\log\left(\rho\right)\leq \overline{\rho\log\left(\rho\right)}, \quad \mbox{a.e. in } (0,1)\times(0,T).
\eeq
Subtracting \eqref{pfsec8.3} by
\eqref{pfsec8.5} and sending $\delta\rightarrow 0$, we have 
\beq\label{pfsec8.7}
\begin{split}
&\int_0^1\overline{\rho\log\left(\rho\right)}(T)-\int_0^1\left(\rho+\frac1j\right)\log\left(\rho+\frac1j\right)(T)\\
=&\int_0^1\rho_0\log\left(\rho_0\right)-\int_0^1\left(\rho_0+\frac1j\right)\log\left(\rho_0+\frac1j\right)\\
&+\int_0^T\int_0^1\rho (u)_x-\lim\limits_{\delta\rightarrow 0}\int_0^T\int_0^1\rho_\delta (u_\delta)_x-\frac1j\int_0^T\int_0^1u_x\log\left(\rho+\frac1j\right).
\end{split}
\eeq
The first two terms of right hand side can be estimated as follows
\beq\label{pfsec8.8}
\begin{split}
&\int_0^T\int_0^1\rho (u)_x-\lim\limits_{\delta\rightarrow 0}\int_0^T\int_0^1\rho_\delta (u_\delta)_x\\
=&\int_0^T\int_0^1\rho (u)_x-\lim\limits_{\delta\rightarrow 0}\int_0^T\int_0^1\rho_\delta \h_\delta^1-\lim\limits_{\delta\rightarrow 0}\int_0^T\int_0^1 A^{-1}_{11}(n_\delta)\rho_\delta^{\gamma+1}\\
=&\int_0^T\int_0^1\rho (u)_x-\int_0^T\int_0^1\rho \h^1-\lim\limits_{\delta\rightarrow 0}\int_0^T\int_0^1 A^{-1}_{11}(n_\delta)\rho_\delta^{\gamma+1}\\
=&\int_0^T\int_0^1\rho A^{-1}_{11}(n)\overline{\rho^{\gamma}}-\lim\limits_{\delta\rightarrow 0}\int_0^T\int_0^1 A^{-1}_{11}(n)\rho_\delta^{\gamma+1}-\lim\limits_{\delta\rightarrow 0}\int_0^T\int_0^1 \left(A^{-1}_{11}(n_\delta)-A^{-1}_{11}(n)\right)\rho_\delta^{\gamma+1}\\
=&\lim\limits_{\delta\rightarrow 0}\int_0^T\int_0^1 A^{-1}_{11}(n)\left(\rho\overline{\rho^{\gamma}}- \rho_\delta^{\gamma+1}\right),
\end{split}
\eeq
where we have used Lemma \ref{lemma5.3} in the second equality, and \eqref{pfsec5.9}, $\gamma>1$, and \eqref{estrho2g} in the last step. Here $\h^1$ is the first element of $\h$, and $A^{-1}_{11}(\cdot)$ is the $(1,1)$ element of inverse matrix $A^{-1}(\cdot)$. By the estimate \eqref{postiveA} and the property of $2\times 2$ matrices,  $A^{-1}_{11}(\cdot)>0$.

Since $\rho, \rho_\delta\geq 0$, it is not hard to verify that 
\beq\notag
(\rho-\rho_\delta)^{\gamma+1}= (\rho-\rho_\delta)^{\gamma}(\rho-\rho_\delta)
\leq \left(\rho^\gamma-\rho_\delta^\gamma\right)(\rho-\rho_\delta).
\eeq
Thus
\beq\label{pfsec8.9}
\begin{split}
&\overline{\lim\limits_{\delta\rightarrow 0}}\int_0^T\int_0^1 A^{-1}_{11}(n)(\rho-\rho_\delta)^{\gamma+1}\\
\leq &\lim\limits_{\delta\rightarrow 0}\int_0^T\int_0^1 A^{-1}_{11}(n) \left(\rho^\gamma-\rho_\delta^\gamma\right)(\rho-\rho_\delta)\\
=&\lim\limits_{\delta\rightarrow 0}\int_0^T\int_0^1 A^{-1}_{11}(n) \left(\rho^{\gamma+1}-\rho^{\gamma}\rho_\delta-\rho_\delta^\gamma\rho+\rho_\delta^{\gamma+1}\right)\\
=&\lim\limits_{\delta\rightarrow 0}\int_0^T\int_0^1 A^{-1}_{11}(n) \left(\rho_\delta^{\gamma+1}-\rho\overline{\rho^{\gamma}}\right)+\lim\limits_{\delta\rightarrow 0}\int_0^T\int_0^1 A^{-1}_{11}(n) \left(\rho^{\gamma+1}-\rho^{\gamma}\rho_\delta-\rho_\delta^\gamma\rho+\rho\overline{\rho^{\gamma}}\right)\\
=&\lim\limits_{\delta\rightarrow 0}\int_0^T\int_0^1 A^{-1}_{11}(n) \left(\rho_\delta^{\gamma+1}-\rho\overline{\rho^{\gamma}}\right).
\end{split}
\eeq
Substituting \eqref{pfsec8.9} into \eqref{pfsec8.8}, we have 
\beq\notag
\int_0^T\int_0^1\rho (u)_x-\lim\limits_{\delta\rightarrow 0}\int_0^T\int_0^1\rho_\delta (u_\delta)_x\leq 0.
\eeq
Combing this inequality with \eqref{pfsec8.7}, we conclude that 
\beq\notag
\begin{split}
&\int_0^1\overline{\rho\log\left(\rho\right)}(T)-\int_0^1\left(\rho+\frac1j\right)\log\left(\rho+\frac1j\right)(T)\\
\leq &\int_0^1\rho_0\log\left(\rho_0\right)-\int_0^1\left(\rho_0+\frac1j\right)\log\left(\rho_0+\frac1j\right)-\frac1j\int_0^T\int_0^1u_x\log\left(\rho+\frac1j\right).
\end{split}
\eeq
Sending $j\rightarrow \infty$, we obtain that 
\beq\notag
\begin{split}
\int_0^1\overline{\rho\log\left(\rho\right)}(T)-\int_0^1\rho\log\left(\rho\right)(T)\leq 0.
\end{split}
\eeq
This and \eqref{pfsec8.6} imply that $\overline{\rho\log\left(\rho\right)}=\rho\log\left(\rho\right)$, combined with \eqref{pfsec8.10},
implies \eqref{pfsec6.2}. 

Combining Lemma \ref{lemma5.4} with Lemma \ref{lemma5.5}, and using the strict convexity of 
$\rho\log\rho$ for $\rho\geq 0$,  we know that 
\beq\notag
\rho_\delta\rightarrow\rho, \quad \mbox{a.e. in } (0,1)\times(0,T).
\eeq
It follows from the Egorov theorem that for any $\epsilon>0$, there is $I_{\epsilon}\subset (0,1)\times(0,T)$ such that $|\big((0,1)\times(0,T)\big)\setminus I_{\epsilon}|<\epsilon$ and 
$$
\sup\limits_{(x,t)\in  I_{\epsilon}}
|\rho_\delta(x,t)-\rho(x,t)|\rightarrow 0.
$$ 
Since $\rho_\delta$ is uniformly bounded in $L^{2\gamma}$, we can estimate
\begin{eqnarray*}
\notag
\int_0^T\int_0^1|\rho_\delta-\rho|^{\gamma}
&\leq& \sup\limits_{(x,t)\in  I_{\epsilon}}
|\rho_\delta(x,t)-\rho(x,t)||I_\epsilon|
+C|\big((0,1)\times(0,T)\big)\setminus I_{\epsilon}|^{\frac12}\|\rho_\delta-\rho\|_{L^{2\gamma}}^{\gamma}
\\
\rightarrow 0, \ \ {\rm{as}}\ \ \delta\rightarrow 0.
\end{eqnarray*}
This implies that $\overline{\rho^\gamma}=\rho^\gamma$ in $(0,1)\times (0,T)$. This completes the proof of Lemma \ref{lemma5.4}.
\endpf



\end{document}